\documentclass[12pt]{article}

\usepackage{amsmath,amsthm}
\usepackage{amssymb}

\def\v{\vert}

\def\b{\ensuremath{\backslash}}
\def\s{\ensuremath{\sigma}}
\def\si{\ensuremath{\sigma}}

\def\u{upstep }

\def\us{upsteps }
\def\ds{downsteps }

\def\un{\underbrace}
\def\t{\textrm}
\def\De{\Delta}

\catcode`\@=11

\thicklines
\newskip\Einheit \Einheit=0.6cm
\newcount\xcoord \newcount\ycoord
\newdimen\xdim \newdimen\ydim \newdimen\PfadD@cke \newdimen\Pfadd@cke
\def\PfadDicke#1{\PfadD@cke#1 \divide\PfadD@cke by2 
\Pfadd@cke\PfadD@cke \multiply\PfadD@cke by2}
\long\def\LOOP#1\REPEAT{\def\BODY{#1}\ITERATE}
\def\ITERATE{\BODY \let\next\ITERATE \else\let\next\relax\fi \next}
\let\REPEAT=\fi
\def\Punkt{\hbox{\raise-2pt\hbox to0pt{\hss\scriptsize$\bullet$\hss}}}

\def\DuennPunkt(#1,#2){\unskip
  \raise#2 \Einheit\hbox to0pt{\hskip#1 \Einheit
          \raise-1.5pt\hbox to0pt{\hss\tiny$\bullet$\hss}\hss}}

\def\NormalPunkt(#1,#2){\unskip
  \raise#2 \Einheit\hbox to0pt{\hskip#1 \Einheit
          \raise-3pt\hbox to0pt{\hss\large$\bullet$\hss}\hss}}
\def\DickPunkt(#1,#2){\unskip
  \raise#2 \Einheit\hbox to0pt{\hskip#1 \Einheit
          \raise-4pt\hbox to0pt{\hss\Large$\bullet$\hss}\hss}}
\def\Kreis(#1,#2){\unskip
  \raise#2 \Einheit\hbox to0pt{\hskip#1 \Einheit
          \raise-4pt\hbox to0pt{\hss\Large$\circ$\hss}\hss}}
\def\Diagonale(#1,#2)#3{\unskip\leavevmode
  \xcoord#1\relax \ycoord#2\relax
      \raise\ycoord \Einheit\hbox to0pt{\hskip\xcoord \Einheit
         \unitlength\Einheit
         \line(1,1){#3}\hss}}
\def\AntiDiagonale(#1,#2)#3{\unskip\leavevmode
  \xcoord#1\relax \ycoord#2\relax \advance\xcoord by -0.05\relax
      \raise\ycoord \Einheit\hbox to0pt{\hskip\xcoord \Einheit
         \unitlength\Einheit
         \line(1,-1){#3}\hss}}
\def\Pfad(#1,#2),#3\endPfad{\unskip\leavevmode
  \xcoord#1 \ycoord#2 \thicklines\ZeichnePfad#3\endPfad\thinlines}
\def\ZeichnePfad#1{\ifx#1\endPfad\let\next\relax
  \else\let\next\ZeichnePfad
    \ifnum#1=1
      \raise\ycoord \Einheit\hbox to0pt{\hskip\xcoord \Einheit
         \vrule height\Pfadd@cke width1 \Einheit depth\Pfadd@cke\hss}%
      \advance\xcoord by 1
    \else\ifnum#1=2
      \raise\ycoord \Einheit\hbox to0pt{\hskip\xcoord \Einheit
        \hbox{\hskip-\PfadD@cke\vrule height1 \Einheit 
width\PfadD@cke depth0pt}\hss}%
      \advance\ycoord by 1
    \else\ifnum#1=3
      \raise\ycoord \Einheit\hbox to0pt{\hskip\xcoord \Einheit
         \unitlength\Einheit
         \line(1,1){1}\hss}
      \advance\xcoord by 1
      \advance\ycoord by 1
    \else\ifnum#1=4
      \raise\ycoord \Einheit\hbox to0pt{\hskip\xcoord \Einheit
         \unitlength\Einheit
         \line(1,-1){1}\hss}
      \advance\xcoord by 1
      \advance\ycoord by -1
    \fi\fi\fi\fi
  \fi\next}
\def\hSSchritt{\leavevmode\raise-.4pt\hbox 
to0pt{\hss.\hss}\hskip.2\Einheit
  \raise-.4pt\hbox to0pt{\hss.\hss}\hskip.2\Einheit
  \raise-.4pt\hbox to0pt{\hss.\hss}\hskip.2\Einheit
  \raise-.4pt\hbox to0pt{\hss.\hss}\hskip.2\Einheit
  \raise-.4pt\hbox to0pt{\hss.\hss}\hskip.2\Einheit}
\def\vSSchritt{\vbox{\baselineskip.2\Einheit\lineskiplimit0pt
\hbox{.}\hbox{.}\hbox{.}\hbox{.}\hbox{.}}}
\def\DSSchritt{\leavevmode\raise-.4pt\hbox to0pt{%
  \hbox to0pt{\hss.\hss}\hskip.2\Einheit
  \raise.2\Einheit\hbox to0pt{\hss.\hss}\hskip.2\Einheit
  \raise.4\Einheit\hbox to0pt{\hss.\hss}\hskip.2\Einheit
  \raise.6\Einheit\hbox to0pt{\hss.\hss}\hskip.2\Einheit
  \raise.8\Einheit\hbox to0pt{\hss.\hss}\hss}}
\def\dSSchritt{\leavevmode\raise-.4pt\hbox to0pt{%
  \hbox to0pt{\hss.\hss}\hskip.2\Einheit
  \raise-.2\Einheit\hbox to0pt{\hss.\hss}\hskip.2\Einheit
  \raise-.4\Einheit\hbox to0pt{\hss.\hss}\hskip.2\Einheit
  \raise-.6\Einheit\hbox to0pt{\hss.\hss}\hskip.2\Einheit
  \raise-.8\Einheit\hbox to0pt{\hss.\hss}\hss}}
\def\SPfad(#1,#2),#3\endSPfad{\unskip\leavevmode
  \xcoord#1 \ycoord#2 \ZeichneSPfad#3\endSPfad}
\def\ZeichneSPfad#1{\ifx#1\endSPfad\let\next\relax
  \else\let\next\ZeichneSPfad
    \ifnum#1=1
      \raise\ycoord \Einheit\hbox to0pt{\hskip\xcoord \Einheit
         \hSSchritt\hss}%
      \advance\xcoord by 1
    \else\ifnum#1=2
      \raise\ycoord \Einheit\hbox to0pt{\hskip\xcoord \Einheit
        \hbox{\hskip-2pt \vSSchritt}\hss}%
      \advance\ycoord by 1
    \else\ifnum#1=3
      \raise\ycoord \Einheit\hbox to0pt{\hskip\xcoord \Einheit
         \DSSchritt\hss}
      \advance\xcoord by 1
      \advance\ycoord by 1
    \else\ifnum#1=4
      \raise\ycoord \Einheit\hbox to0pt{\hskip\xcoord \Einheit
         \dSSchritt\hss}
      \advance\xcoord by 1
      \advance\ycoord by -1
    \fi\fi\fi\fi
  \fi\next}
\def\Koordinatenachsen(#1,#2){\unskip
 \hbox to0pt{\hskip-.5pt\vrule height#2 \Einheit width.5pt depth1 
\Einheit}%
 \hbox to0pt{\hskip-1 \Einheit \xcoord#1 \advance\xcoord by1
    \vrule height0.25pt width\xcoord \Einheit depth0.25pt\hss}}
\def\Koordinatenachsen(#1,#2)(#3,#4){\unskip
 \hbox to0pt{\hskip-.5pt \ycoord-#4 \advance\ycoord by1
    \vrule height#2 \Einheit width.5pt depth\ycoord \Einheit}%
 \hbox to0pt{\hskip-1 \Einheit \hskip#3\Einheit 
    \xcoord#1 \advance\xcoord by1 \advance\xcoord by-#3 
    \vrule height0.25pt width\xcoord \Einheit depth0.25pt\hss}}
\def\Gitter(#1,#2){\unskip \xcoord0 \ycoord0 \leavevmode
  \LOOP\ifnum\ycoord<#2
    \loop\ifnum\xcoord<#1
      \raise\ycoord \Einheit\hbox to0pt{\hskip\xcoord 
\Einheit\Punkt\hss}%
      \advance\xcoord by1
    \repeat
    \xcoord0
    \advance\ycoord by1
  \REPEAT}
\def\Gitter(#1,#2)(#3,#4){\unskip \xcoord#3 \ycoord#4 \leavevmode
  \LOOP\ifnum\ycoord<#2
    \loop\ifnum\xcoord<#1
      \raise\ycoord \Einheit\hbox to0pt{\hskip\xcoord 
\Einheit\Punkt\hss}%
      \advance\xcoord by1
    \repeat
    \xcoord#3
    \advance\ycoord by1
  \REPEAT}
\def\Label#1#2(#3,#4){\unskip \xdim#3 \Einheit \ydim#4 \Einheit
  \def\lo{\advance\xdim by-.5 \Einheit \advance\ydim by.5 \Einheit}%
  \def\llo{\advance\xdim by-.25cm \advance\ydim by.5 \Einheit}%
  \def\loo{\advance\xdim by-.5 \Einheit \advance\ydim by.25cm}%
  \def\o{\advance\ydim by.25cm}%
  \def\ro{\advance\xdim by.5 \Einheit \advance\ydim by.5 \Einheit}%
  \def\rro{\advance\xdim by.25cm \advance\ydim by.5 \Einheit}%
  \def\roo{\advance\xdim by.5 \Einheit \advance\ydim by.25cm}%
  \def\l{\advance\xdim by-.30cm}%
  \def\r{\advance\xdim by.30cm}%
  \def\lu{\advance\xdim by-.5 \Einheit \advance\ydim by-.6 \Einheit}%
  \def\llu{\advance\xdim by-.25cm \advance\ydim by-.6 \Einheit}%
  \def\luu{\advance\xdim by-.5 \Einheit \advance\ydim by-.30cm}%
  \def\u{\advance\ydim by-.30cm}%
  \def\ru{\advance\xdim by.5 \Einheit \advance\ydim by-.6 \Einheit}%
  \def\rru{\advance\xdim by.25cm \advance\ydim by-.6 \Einheit}%
  \def\ruu{\advance\xdim by.5 \Einheit \advance\ydim by-.30cm}%
  #1\raise\ydim\hbox to0pt{\hskip\xdim
     \vbox to0pt{\vss\hbox to0pt{\hss$#2$\hss}\vss}\hss}%
}
\catcode`\@=12

\begin{document}
\theoremstyle{plain}
\newtheorem{theorem}{Theorem}
\newtheorem{lemma}[theorem]{Lemma}
\newtheorem{prop}[theorem]{Proposition}
\newtheorem{cor}[theorem]{Corollary}

\begin{center}
{\Large
  A recursive bijective approach to counting permutations containing 
  3-letter patterns                         \\ 
}
\vspace{10mm}
DAVID CALLAN  \\
Department of Statistics  \\
University of Wisconsin-Madison  \\
1210 W. Dayton St   \\
Madison, WI \ 53706-1693  \\
{\bf callan@stat.wisc.edu}  \\
\vspace{5mm}
\today
\end{center}

\vspace{5mm}

We present a method, illustrated by several examples, to find explicit 
counts of permutations containing a given multiset of three letter 
patterns. The method is recursive, depending on bijections to reduce 
to the case of a smaller multiset, and involves a consideration of separate 
cases according to how
the patterns overlap. Specifically, we use the method (i) to 
provide combinatorial proofs of  B\'{o}na's formula $\binom{2n-3}{n-3}$ 
for the number of $n$-permutations containing one 132 pattern and Noonan's formula
$\frac{3}{n}\binom{2n}{n+3}$ for one 123 pattern, 
(ii) to express the number of $n$-permutations containing exactly $k$ 123 
patterns in terms of ballot numbers for $k\le 4$, and (iii) 
to express the number of 123-avoiding $n$-permutations containing exactly $k$ 
132 patterns as a linear combination of powers of 2, also for $k\le 4$. The 
results strengthen the conjecture that the counts are algebraic for 
all $k$.

\vspace*{10mm}
{\Large \textbf{1 \quad Introduction}  }

In the context of pattern-containing permutations, a 
\emph{permutation} $\pi$ is a list (or \emph{word}) 
$\pi_{1}\,\pi_{2}\,\ldots\,\pi_{n}$ of distinct 
positive integers called letters. Viewed as a bijection, $\pi$ 
sends its $i$th smallest letter to $\pi_{i}$. The \emph{reverse} of 
a permutation  $\pi_{1}\,\pi_{2}\,\ldots\,\pi_{n}$ is 
$\pi_{n}\,\pi_{n-1}\,\ldots\,\pi_{1}$ and, when the letters of $\pi$ 
are $\{1,2,\ldots,n\}$, its \emph{complement} is 
$n+1-\pi$ (termwise).
The \emph{reduced form} of a permutation $\pi$, denoted reduce($\pi$), 
is obtained by replacing its smallest element by 1, its next smallest 
by 2, and so on. For example, reduce$(9\,8\,2\,4\,6) =5\,4\,1\,2\,3$. We use 
$[n]$ to denote $\{1,2,\ldots,n\}$ and $[m,n]$ to denote the interval 
of integers $\{m,m+1,\ldots,n\}$.
A $k$-letter \emph{pattern} is simply a permutation
$\tau$ on $[k]$. An occurrence of $\tau$ in a permutation $\pi$ is 
a (scattered) subword $\pi_{i_{1}}\,\pi_{i_{2}}\,\ldots\pi_{i_{k}}$ of $\pi$ whose 
reduced form is $\tau$. For example, $4\,3\,1\,2$ has two occurrences of the 
pattern 321, namely $4\,3\,1$ and $4\,3\,2$, but it is 132-avoiding.
We will refer to the actual letters in an occurrence of a 321 pattern as $c,\,b,\,a$ 
respectively. Similarly for 132 patterns.

Much work has recently been done counting permutations that contain 
a specified number of 
occurrences of one or more patterns.  See \cite{wilf02} for a general survey 
and \cite{fulmek02} for results on 3 letter patterns, the subject of 
this paper.
We let $P_{\tau}^{(i)}(n)$ denote the set of 
permutations on $[n]$ containing exactly $i$ occurrences of the 
pattern $\tau$. For example, $P_{21}^{(1)}(3)=\{1\,3\,2,\ 2\,1\,3\}$.
Of course, $\v P^{(i)}_{123} \v = \v P^{(i)}_{321} \v$ (reversing 
permutations carries patterns of one kind to the other).
In section 2, we develop results counting various classes of lattice 
paths. These results count, via bijection, classes of 
321-avoiding permutations that arise in subsequent sections.

\vspace*{10mm}

{\Large \textbf{2 \quad Lattice Paths}  }

Our starting point is a known formula for counting first quadrant lattice 
paths of upsteps and downsteps.
A \emph{path} is a sequence of contiguous upsteps $(1,1)$ and 
downsteps $(1,-1)$ 
that starts at the origin. A \emph{first quadrant} path is one that 
stays (weakly) in the first quadrant. 
A \emph{balanced} path is one 
with an equal number of upsteps and downsteps, $n$, the \emph{semilength} of 
the path. A \emph{Dyck} path is a balanced first quadrant path. 
Illustrated are a first quadrant path of 6 upsteps and 4 downsteps 
and a Dyck path of semilength 4.

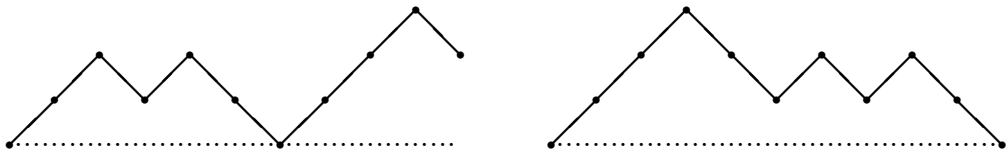
\[
\Pfad(-11,-1),3343443334\endPfad
\SPfad(-11,-1),1111111111\endSPfad
\DuennPunkt(-11,-1)
\DuennPunkt(-10,0)
\DuennPunkt(-9,1)
\DuennPunkt(-8,0)
\DuennPunkt(-7,1)
\DuennPunkt(-6,0)
\DuennPunkt(-5,-1)
\DuennPunkt(-4,0)
\DuennPunkt(-3,1)
\DuennPunkt(-2,2)
\DuennPunkt(-1,1)
\Pfad(1,-1),3334434344\endPfad
\SPfad(1,-1),1111111111\endSPfad
\DuennPunkt(1,-1)
\DuennPunkt(2,0)
\DuennPunkt(3,1)
\DuennPunkt(4,2)
\DuennPunkt(5,1)
\DuennPunkt(6,0)
\DuennPunkt(7,1)
\DuennPunkt(8,0)
\DuennPunkt(9,1)
\DuennPunkt(10,0)
\DuennPunkt(11,-1)
\]

Both have \emph{height} 3 
and one return to the $x$-axis. The Dyck path's only return is at the 
end; it has no \emph{interior} returns. The first path has 3 ascents 
(lengths of maximal sequences of contiguous upsteps)   
2,1,3 and 3 descents 1,2,1. We say its 
\emph{ascent sequence} is $(2,1,3)$ and its \emph{descent sequence} is $(1,2,1)$. 
The ascent and descent sequences for the Dyck path shown are (3,1,1) and (2,1,2) 
respectively.
Clearly, the requirements to be a valid ascent-descent sequence pair 
for a Dyck path are (i) positive integer entries, (ii) same length, 
(iii) same sum, and (iv) each partial sum of the first weakly exceeds 
the corresponding partial sum of the second.
The number of Dyck $n$-paths (first quadrant balanced paths of semilength $n$) 
is well known to be the Catalan 
number $C_{n}$ and general first quadrant paths are counted by 
generalized Catalan numbers (ballot numbers) as follows.

Let $C(x)=\sum_{n=0}^{\infty}C_{n}x^{n}$ denote the (ordinary) generating 
function for the Catalan numbers and set 
$C^{(k)}_{n}=[x^{n}]C(x)^{k}$ (the ballot numbers). Thus 
$C_{n}=C^{(1)}_{n}$. Since 
multiplication of generating functions corresponds to convolution of 
their coefficient sequences, $C^{(k)}_{n}$ is the number of 
sequences of length $k$ of (possibly empty) Dyck paths 
whose total semilength is $n$.
By joining the $k$ paths together with intervening upsteps---a 
reversible procedure---we see that $C^{(k)}_{n}$ is the number of 
first quadrant paths of $n+k-1$ \us and $n$ downsteps. The marvelous 
Andr\'{e} reflection principle, illustrated below, can be used to 
count these paths.
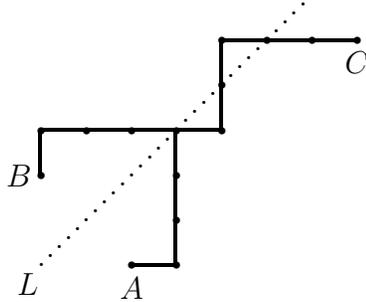
\[
\thinlines 
\SPfad(-4,0),333333\endSPfad
\Pfad(-4,2),2111122111\endPfad
\Pfad(-2,0),1222122111\endPfad
\DuennPunkt(-4,2)
\DuennPunkt(-4,3)
\DuennPunkt(-3,3)
\DuennPunkt(-2,3)
\DuennPunkt(-2,0)
\DuennPunkt(-1,0)
\DuennPunkt(-1,1)
\DuennPunkt(-1,2)
\DuennPunkt(-1,3)
\DuennPunkt(0,3)
\DuennPunkt(0,4)
\DuennPunkt(0,5)
\DuennPunkt(1,5)
\DuennPunkt(2,5)
\DuennPunkt(3,5)
\Label\u{A}(-2,0)
\Label\l{B}(-4,2)
\Label\lu{L}(-3.8,0.2)
\Label\u{C}(3,5)
\]

\noindent The number of paths from $A$ to $C$ that meet line $L$ is the same as 
the total number of paths from $B$ to $C$ (reflect across the first 
point of contact with $L$). Hence, the number of paths from $A$ to $C$ 
that avoid $L$ is the total number from $A$ to $C$ minus the total number from $B$ to $C$.
This yields 
$C^{(k)}_{n}=\binom{2n+k-1}{n}-\binom{2n+k-1}{n-1}=\frac{k}{2n+k}\binom{2n+k}{n}$.
Another nice 
bijective proof that the number of first quadrant paths of $n+k-1$ \us and $n$ downsteps is 
$\frac{k}{2n+k}\binom{2n+k}{n}$ is given in \cite{woan}. In brief, 
consider the $k\binom{2n+k}{n}$ paths of $n+k$ upsteps and $n$ 
downsteps in which one of the $k$ high points has been marked. (The 
high points are the leftmost points on the path at height  
$h,\,h-1,\,\ldots, h-k+1$ respectively, where $h$ is the height of the 
path.) Then the parameter $\nu=$ $x$-coordinate 
of the marked high point is uniformly distributed over its possible 
values $1,2,\ldots,2n+k$. The kicker is that the paths with 
$\nu=2n+k$ are precisely the paths that do not reach their highest 
level until the last upstep; deleting this step and rotating 
$180^{\circ}$, 
they are first quadrant paths of $n+k-1$ \us and $n$ downsteps. So we 
have $C^{(k)}_{n}=\frac{k}{2n+k}\binom{2n+k}{n}$. 

We will make frequent use of the convolution identity 
$\sum_{k=0}^{n} C^{(r)}_{k}C^{(s)}_{n-k}=C^{(r+s)}_{n}$. From the 
above interpretation of $C^{(k)}_{n}$, we have the following simple 
results.
\begin{lemma}
	\label{01}
	The number of 
    first quadrant paths consisting of $m$ upsteps and $n$ downsteps is 
    $C^{(m-n+1)}_{n}$ where 
    $C^{(k)}_{n}=\binom{2n+k-1}{n}-\binom{2n+k-1}{n-1}=\frac{k}{2n+k}\binom{2n+k}{n}$.
\end{lemma}
\begin{lemma}
	$\sum_{j=0}^{\min(m,k-m-1)}C^{(k-2j)}_{j}=\binom{k-1}{m}\qquad m\ge 
	0,\,k\ge 1$
	\label{basic}
\end{lemma}
\begin{proof}
	Using $C^{(k)}_{n}=\binom{2n+k-1}{n}-\binom{2n+k-1}{n-1}$, the left 
	side is a telescoping sum.
\end{proof}
\begin{lemma}
	\[
   \begin{array}{rrcll}
	(i) &  C^{(k)}_{n} - C^{(k-1)}_{n} & =  & C^{(k+1)}_{n-1} &\qquad \ n,k \ge 1 \\
	(ii) &  C^{(k)}_{n}-C^{(k)}_{n-1} & = & C^{(k-2)}_{n}+C^{(k+1)}_{n-1}  & 
	\qquad n \ge 1,\ k \ge 2 \\
	(iii) &  \sum_{j=0}^{n}C^{(k+j)}_{n-j} &  = &  C^{(k+1)}_{n} & \qquad \ 
	n,k \ge 0
   \end{array}
   \]
   \label{05}
\end{lemma}
\emph{Proof.}
Counting first quadrant paths by last step (up or down) yields 
$C^{(k)}_{n}=C^{(k-1)}_{n} + 
C^{(k+1)}_{n-1}$, and (i) follows.
Repeating this decomposition on the first summand
yields (ii). If instead we decompose the second summand and iterate, we get 
(iii) (after replacing $k$ by $k+1$).

A direct combinatorial proof of (iii) can also be given. 
Count the $C^{(k+1)}_{n}$ first quadrant paths of 
$n+k$ \us and $n$ \ds by number of returns to the $x$-axis. Given 
such a path with $j \in [0,n]$ returns, delete the initial \u and 
all $j$ \ds that return the path to the $x$-axis. For example, with 
$j=k=2$ ($D_{1},D_{2}$ denote Dyck paths and $P$ denotes a first 
quadrant path with one more \u than downstep),
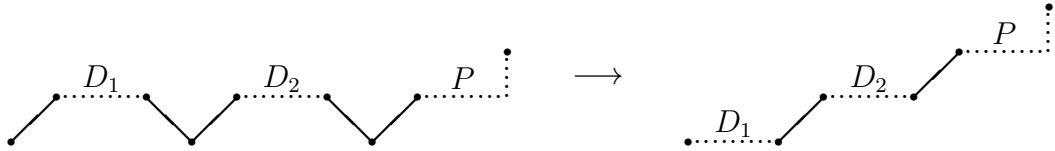
\[
\Pfad(-12,0),3\endPfad
\SPfad(-11,1),11\endSPfad
\SPfad(-7,1),11\endSPfad
\SPfad(-3,1),112\endSPfad
\SPfad(3,0),11\endSPfad
\SPfad(6,1),11\endSPfad
\SPfad(9,2),112\endSPfad
\Pfad(-9,1),43\endPfad
\Pfad(-5,1),43\endPfad
\Pfad(5,0),3\endPfad
\Pfad(8,1),3\endPfad
\DuennPunkt(-12,0)
\DuennPunkt(-11,1)
\DuennPunkt(-9,1)
\DuennPunkt(-8,0)
\DuennPunkt(-7,1)
\DuennPunkt(-5,1)
\DuennPunkt(-4,0)
\DuennPunkt(-3,1)
\DuennPunkt(-1,2)
\DuennPunkt(3,0)
\DuennPunkt(5,0)
\DuennPunkt(6,1)
\DuennPunkt(8,1)
\DuennPunkt(9,2)
\DuennPunkt(11,3)
\Label\o{\longrightarrow}(1,1)
\Label\o{D_{1}}(-10,1)
\Label\o{D_{2}}(-6,1)
\Label\o{P}(-2,1)
\Label\o{D_{1}}(4,0)
\Label\o{D_{2}}(7,1)
\Label\o{P}(10,2)
\] 
This gives a bijection to the first quadrant paths consisting of 
$n+k-1$ \us and $n-j$ downsteps, counted by $C^{(k+j)}_{n-j}$. The 
inverse bijection is : insert an upstep at the start and (if $j>0$) 
locate the rightmost upstep at levels 
$1,2,\ldots,j$ and insert a downstep just before each of these upsteps.

A short 
digression: identity (iii) can be expressed as a property of the 
so-called Catalan triangle $\left(C^{(k)}_{n-k}\right)_{n,k\ge 0}$:
\[
\bordermatrix{
n\backslash k   & 0  & 1   & 2    &  3   &  4   &  5   &  6   &   7  & 8    \cr
0  & 1  &     &      &      &      &      &      &      &      \cr
1  & 0  & 1   &      &      &      &      &      &      &      \cr
2  & 0  & 1   & 1    &      &      &      &      &      &      \cr
3  & 0  & 2   & 2    &  1   &      &      &      &      &      \cr
4  & 0  & 5   & 5    &  3   &   1  &      &      &      &      \cr
5  & 0  & 14  &14    &  9   &   4  & 1    &      &      &      \cr
6  & 0  & 42  &42    & 28   &  14  & 5    & 1    &      &      \cr
7  & 0  & 132 &132   & 90   &  48  & 20   & 6    & 1    &      \cr
8  & 0  & 429 &429   & 297  & 165  & 75   & 27   & 7    &  1   \cr 
}
\]
Pick any column. Sum the entries of each row starting in that column. 
The results will be the entries in the next column. The inverse of this matrix is 
$\left( (-1)^{n-k}\binom{k}{n-k}\right)_{n,k\ge 0}$ and in fact the nonzero 
entries in row $n$ are just the nonzero coefficients
of $q_{n}(x)=x^{\frac{n}{2}}U_{n}(\frac{1}{2\sqrt{x}})$ where $U_{n}(x)$ is 
the Chebychev polynomial of the second kind defined by $\frac{\sin( 
n+1)\theta}{\sin\theta }=U_{n}(\cos \theta)$. This polynomial 
$q_{n}(x)$ crops up as the generating function of Dyck paths of 
bounded height.

\textbf{Fact \cite[Theorem 2]{kratt} }\quad 
	Let $u(n,h)$ denote the number of Dyck $n$-paths of height $\le h$. 
	Then
	\[
	\sum_{n\ge 0}^{}u(n,h)x^{n}=\frac{q_{h}(x)}{q_{h+1}(x)}.
	\]
 
This result can be generalized to count paths bounded by two given 
horizontal lines and ending at a given height.

\textbf{Fact}\quad 
	Let $v(n,h,r,s)$ denote the number of paths of $n+h-1$ upsteps and 
	$n$ downsteps that start at the origin and are weakly bounded by the 
	line $r$ units below the initial point $($i.e. the line $y=-r)$ and the line $s$ 
	units above the terminal point $($i.e. the line $y=s+h-1)$. Then
	\[
		\sum_{n\ge 0}^{}v(n,h,r,s)x^{n}=\frac{q_{r}(x)q_{s}(x)}{q_{r+s+h}(x)}.
	\]
\newpage  
\noindent Two special cases: 

(i) $r=s=0$. The generating function for paths that 
are bounded by two horizontal lines $h-1$ units apart, that start on 
the lower line and end on the upper line is $\frac{1}{q_{h}(x)}$.

(ii) $ r=h-1,\,s=0$. The generating function for paths that 
are bounded by two horizontal lines $2(h-1)$ units apart, that start midway 
between the lines and end on the upper line is 
$\frac{q_{h-1}(x)}{q_{2h-1}(x)}=\frac{1}{p_{h}(x)}$ where 
$p_{h}(x)=2x^{\frac{n}{2}}T_{h}(\frac{1}{2\sqrt{x}})$ and $T_{n}$, 
defined by $T_{n}(\cos \theta)=\cos n\theta$, is the Chebychev 
polynomial of the first kind. 

Thus both kinds of Chebychev polynomials 
have lattice path interpretations. End of digression.

The next few results count various classes of Dyck paths to serve as 
the induction basis for recursively counting pattern-containing permutations.
\begin{lemma}
	\label{02}
	The number of Dyck $n$-paths with first ascent $=k$ is 
	$C^{(k)}_{n-k}$ and the number with first ascent $\ge k$ is 
	$C^{(k+1)}_{n-k}$.
\end{lemma}
\begin{proof}
	A path with first ascent $=k$ begins with $k$ upsteps followed by a 
	downstep. Deleting these $k$ upsteps and 1 downstep and reading the resulting path 
	backwards gives a bijection to first quadrant paths of $n-1$ upsteps 
	and $n-k$ downsteps, counted by $C^{(k)}_{n-k}$. For paths with first ascent 
	$\ge k$, delete the first $k$ upsteps and read backwards to get a bijection 
    to first quadrant paths of $n$ upsteps 
	and $n-k$ downsteps, counted by $C^{(k+1)}_{n-k}$. 
\end{proof}
\begin{lemma}
	\label{03}
	$($i\,$)$ The number of Dyck $n$-paths with first ascent $\ge r$, at least one 
	interior return, and last descent $\ge s$ is $C^{(r+s+1)}_{n-r-s}$.\\
	$($ii\,$)$ The total number of Dyck $n$-paths with first ascent $\ge r$ 
	and last descent $\ge s$ is 
	$\sum_{j=0}^{\min(r,s)}C^{(r+s+1-2j)}_{n-r-s+j}$
\end{lemma}
\begin{proof}
	(i) Given such a path, remove the first downstep that returns the path 
	to the $x$-axis and transfer the first upstep to the vacated location, then place the removed 
	downstep at the end of the path. Repeat this process $r$ times 
	altogether. This gives a bijection to Dyck $n$-paths whose last 
	descent is $\ge r+s$ and, reading the paths backwards, to Dyck $n$-paths 
	whose first ascent is $\ge r+s$, counted by  $C^{(r+s+1)}_{n-r-s}$ 
	(Lemma \ref{02}). \\
	(ii) Using (i), count by the lowest point on the path.
\end{proof}
\begin{lemma}
	For $a,d,i\ge 1$ and $n\ge i+d$, the Dyck $n$-paths  with first 
	ascent $a$, last descent 
	$d$, and first $i$ \emph{nonfinal} descents 
	all $=1$ are equinumerous with the Dyck $n$-paths with first ascent 
    $a+i$ and last descent $d$.
	\label{nonfinal}
\end{lemma}
\begin{proof}
	There is a  cut-and-paste bijection from the first set to the second:
	simply transfer the upstep immediately following each of the first $i$ 
    downsteps to the start of the path. It's inverse is: remove $i$ 
	upsteps from the start of the path, and insert them, one apiece, 
	immediately after each of the first $i$ downsteps (which will not be 
	steps of the final descent because $n\ge i+d$).
\end{proof}
\begin{lemma}
    The Dyck $n$-paths  with first ascent $\ge r$ and first $s$ nonfinal descents 
	all $=1$ are equinumerous with the Dyck $n$-paths with first ascent 
    $\ge r+s$, and hence counted by  $C^{(r+s+1)}_{n-r-s}$.
	\label{oddball}
\end{lemma}
\begin{proof}
	Sum the previous lemma's count over $a\ge r$ and $d\ge 1$ with $i:=s$.
\end{proof}

\begin{lemma}
	The number of Dyck $n$-paths  with first ascent $\ge r$ and last $s-1$ 
	noninitial ascents all $=1$ is
	\[
	\sum_{j=0}^{m}C^{(r+s+1-2j)}_{n-r-s+j}-
	\sum_{j=2}^{m}\binom{r+s-2j}{r-j}C^{(0)}_{n-r-s+j}
	\]
	where $m=\min(r,s)$.
	\label{oddball1}
\end{lemma}
\begin{proof}
    Reading the paths backwards, they correspond to Dyck $n$-paths 
    with first $s-1$ nonfinal descents all $=1$ and last descent $\ge 
    r$. For $n\ge r+s-1$, the latter correspond by Lemma \ref{nonfinal}  
	to Dyck $n$-paths  with first ascent $\ge s$ and last descent $\ge 
	r$ and hence  are counted by the first sum (Lemma \ref{03}\,(ii)). The second sum 
    is a correction for the case $n<r+s-1$ (using Lemma \ref{basic}).
\end{proof}
There is another formulation of the previous result that avoids subtraction:
\begin{lemma}
	The number of Dyck $n$-paths  with first ascent $\ge r$ and last $s-1$ 
	noninitial ascents all $=1$ is
	\[
      C^{(r+s)}_{n+1-r-s}+\sum_{k=0}^{r+s-4}\
	  \sum_{j=0}^{r+s-4-k}\binom{k}{r-2-j}C^{(r+s-2-k)}_{n-r-s+1}.
	\]
\end{lemma}
\begin{proof}
	Let $F(r,s)=\sum_{n\ge 
	0}^{}\Big(C^{(r+s)}_{n+1-r-s}+\sum_{k=0}^{r+s-4}\ \sum_{j=0}^{r+s-4-k}
	\binom{k}{r-2-j}C^{(r+s-2-k)}_{n+1-r-s}\Big)x^{n}$ and 
	$G(r,s)=
	\sum_{n\ge 0}^{}\Big(\sum_{j=0}^{\min(r,s)}C^{(r+s+1-2j)}_{n-r-s+j}-
	\sum_{j=2}^{\min(r,s)}\binom{r+s-2j}{r-j}C^{(0)}_{n-r-s+j}\Big)x^{n}$ denote 
	the generating functions for the two expressions whose equality we 
	wish to establish. We will show that $F$ and $G$ both satisfy the 
	recurrence
	\begin{equation}
		u(r,s)=xu(r-1,s)+xu(r,s-1)+x^{r+s-1}C(x)^{r+s-2}\big(1-C(x)\big)^{2}\qquad 
		r,s\ge 2
		\label{eq:1}
	\end{equation}
	with initial condition
	\[
	u(r,1)=u(1,r)=x^{r}C(x)^{r+1}\qquad r \ge 1
	\]
	(and hence $F(r,s)=G(r,s)$ for all $r,s\ge 1$).
	
	From the definition of $C^{(m)}_{n}$, we have 
	\[
	F(r,s)=x^{r+s-1}C(x)^{r+s}+\sum_{k=0}^{r+s-4}\ \sum_{j=0}^{r+s-4-k}
	\binom{k}{r-2-j}x^{r+s-1}C(x)^{r+s-2-k}
	\]
	and 
	\[
	G(r,s)=\sum_{j=0}^{\min(r,s)} 
	x^{r+s-j}C(x)^{r+s+1-2j}-\sum_{j=2}^{\min(r,s)}\binom{r+s-2j}{r-j} 
	x^{r+s-j}.
	\]
	It is routine to verify $G(r,1)=x^{r}C(x)^{r+1}$ using 
	the basic identity $xC(x)^{2}=C(x)-1$ for the Catalan numbers 
	generating function. Since both $G$ and the recurrence are symmetric 
	in $r,s$ it suffices to verify that $G$ satisfies the recurrence 
	for $r \ge s\ (\ge 2)$. First suppose $r>s$; this guarantees 
	$\min(r-1,s)=s$. (With a little care, detailed later, the same 
	calculation works in case $r=s$.) Adjusting summations to terminate 
	at $s$ and using the basic identity 
	$\binom{a}{b-1}+\binom{a}{b}=\binom{a+1}{b}$, integer $a,b \ge 0$, 
	the right hand side of (\ref{eq:1}) with $u=G$ can be written
	\[
	2\sum_{j=0}^{s}x^{r+j}C(x)^{r-s+2j}-
	x^{r}C(x)^{r-s}-\sum_{j=2}^{s}\binom{r+s-2j}{r-j}x^{r+s-j}+
	x^{r+s-1}C(x)^{r+s-2}\big(1-C(x)\big)^{2}.
	\]
	The left hand side is $G(r,s)$, and the difference is 
	\[
	x^{r}C(x)^{r-s}\big(2-C(x)\big)
	\sum_{j=0}^{s}\big(xC(x)^{2}\big)^{j}-x^{r}C(x)^{r-s}+
	x^{r+s-1}C(x)^{r+s-2}\big(1-C(x)\big)^{2}.
	\]
	Evaluating the sum and using $2-C(x)=1-xC(x)^{2}$, this difference 
	simplifies to 
	$x^{r+s-1}C(x)^{r+s-2}\Big(\big(1-C(x)\big)^{2}-x^{2}C(x)^{4}\Big)=0$.
	
	In case $r=s$, the two sums in $G(r-1,s)$ terminate at $s-1$ rather 
	than $s$, but the contributions from $j=s$ cancel provided we 
	define $\binom{-1}{-1}=1$. The adjustment of summations to end at $s$ 
	holds provided $\binom{-1}{0}=0$. Making these definitions, the 
	identity $\binom{a}{b-1}+\binom{a}{b}=\binom{a+1}{b}$ holds also for 
	$a=-1,\,b=0$ and the entire calculation above goes through for $r=s$.
	
	As for $F$, reversing the inner summation and then replacing the 
	lower binomial coefficient parameter by its complement shows that $F$ is symmetric in $r$ 
	and $s$. We have  
	$F(r,1)=x^{r}C(x)^{r+1}$ because the double sum makes no 
	contribution when $s=1$. The right hand side of (\ref{eq:1}) with 
	$u=F$ is, after combining terms,
	 \begin{eqnarray*}
	 	 &  &  	x^{r+s-1}\Big(C(x)^{r+s}+C(x)^{r+s-2}+
	\sum_{\substack{j,k\ge 0, \\ j+k \le 
	r+s-5}}\binom{k+1}{r-2-j}C(x)^{r+s-k-3}\Big) \\
	 	 & = &  	x^{r+s-1}\Big(C(x)^{r+s}+C(x)^{r+s-2}+
	\sum_{\substack{j\ge 0,K\ge 1, \\ j+K \le 
	r+s-4}}\binom{K}{r-2-j}C(x)^{r+s-K-2}\Big)  \\
	 	 & = &  x^{r+s-1}\Big(C(x)^{r+s}+
	\sum_{\substack{j,K\ge 0 \\ j+K \le 
	r+s-4}}\binom{K}{r-2-j}C(x)^{r+s-K-2}\Big) \\
	 	 & = & F(r,s)\quad \textrm{=\quad the left hand side of (\ref{eq:1}). }
	 \end{eqnarray*}
	 (The hypothesis $r,s\ge 2$ is needed for the next to last equality.)	
	
\end{proof}
\begin{lemma}
	\label{13}
	\[
	\sum_{k=0}^{m-1}\binom{2m-2k}{m-1-k}C_{k}=\binom{2m+1}{m-1}.
	\]
\end{lemma}
\begin{proof}
	The right side counts lattice paths of $m+2$ \us and $m-1$ downsteps. Given 
	such a path $P$, note that it ends at height 3 and split it into two 
	subpaths $P_{1},P_{2}$ by deleting the last upstep that carries 
	$P$ to height 3. Then $P_{1}$ is a path with 2 more upsteps  than downsteps  
	and $P_{2}$ is a Dyck path. Counting by the semilength $k$ of this 
	Dyck path yields the claimed identity.
\end{proof}

\newpage 

{\Large \textbf{3 \quad 321-avoiding permutations}  }

Several bijections are known  \cite{kratt,astrid} among  
$\,P^{(0)}_{321}(n),\,P^{(0)}_{132}(n)$ and Dyck paths of semilength 
$n$, proving that  $\v 
P^{(0)}_{321}(n) \v =\v P^{(0)}_{132}(n) \v = C_{n}$. 
We need to examine one of them. Define the difference operator $\De$ on a 
finite sequence, as usual, by 
$\De\big((u_{i})_{i=1}^{n}\big)=(u_{i}-u_{i-1})_{i=2}^{n}$. 
\begin{theorem}[Krattenthaler\cite{kratt}]
	\label{06}
	There is a bijection $\pi \to P$ from $P^{(0)}_{321}(n)$ to Dyck 
	$n$-paths that sends $\pi_{1}$ to the first ascent of $P$.
\end{theorem}
The bijection can be described as follows (see \cite{fulmek02} for a nice 
pictorial description). Given $\pi \in 
P^{(0)}_{321}(n)$, let $1=i_{1}<i_{2}<\ldots<i_{k}$ denote the 
positions (locations) of the successive record highs (left to right 
maxima) in $\pi$. Then $\De(0,\pi_{i_{1}},\pi_{i_{2}},\ldots,\pi_{i_{k}}=n)$ is the 
ascent sequence of $P$ and $\De(i_{1}=1,i_{2},\ldots,i_{k},n+1)$ is 
the descent sequence of $P$. For example, 
$2\,1\,4\,7\,3\,5\,6\in P^{(0)}_{321}(7)$ 
has record highs of 2,4,7 in positions 1,3,4 yielding ascent 
sequence $\De(0,2,4,7)=(2,2,3)$ and descent sequence 
$\De(1,3,4,8)=(2,1,4)$. Since $i_{1}=1$, we see that $\pi_{1}$ is the 
first ascent of $P$.

We use this bijection to count certain classes of 321-avoiding 
permutations.
\begin{lemma}
	\label{07}
	$\#\{\pi \in P^{(0)}_{321}(n):\ \pi_{1}=k\} = C^{(k)}_{n-k}$ and \ 
	$\#\{\pi \in P^{(0)}_{321}(n):\ \pi_{1}\ge k\} = C^{(k+1)}_{n-k}$
\end{lemma}
\begin{proof}
	Immediate from the bijection and Lemma \ref{02}.
\end{proof}
\begin{lemma}
	\label{08}
	For $m \in [n]$, the $321$-avoiding permutations on $[n]$  with the following 
	properties are equinumerous and all are counted by 
	$C^{(m+1)}_{n-m}$. 
	\begin{enumerate}
		\item $1$ occurs no earlier than position $m$
		\item first entry is $\ge m$
		\item $n$ occurs no later than the $m$th position from the end
		\item last entry is $\le n+1-m$
	\end{enumerate}
\end{lemma}
\begin{proof}
	The map ``reverse and 
	complement'' is a bijection from the first set to the third set 
	and also from the second to the fourth.
	For a permutation $\pi$ in the first set, the first $m-1$ entries 
	are increasing and the first entry is $\ge 2$. So the table of 
	record highs and their locations begins 
	$\begin{smallmatrix}
	\t{location:} & 1 & 2 & \ldots & m-1 \\
	\t{record high:} & \pi_{1} & \pi_{2} & \ldots & \pi_{m-1}
	\end{smallmatrix}$
	with $\pi_{1}\ge 2$. Under Krattenthaler's bijection, the descent 
	sequence of the corresponding Dyck path, 
	$\De(1,2,\ldots,m-1,\,other\ locations)$, begins $\un{1,1,\ldots,1}_{m-2}$ 
	and its ascent sequence, $\De(0,\pi_{1},\ldots)$, begins with 
	$\pi_{1} \ge 2$. This means we get all Dyck paths with first ascent 
	$\ge 2$ and the first $m-2$ descents all $=1$. The second set is sent 
	to Dyck $n$-paths with first ascent $\ge m$.  Lemma \ref{oddball} now 
	yields the result.
\end{proof}
\begin{lemma}
	\label{08a}
	Let $n>i\ge 1,\ k \in [n]$ and consider the $321$-avoiding permutations on $[n]$ with first entry 
	$k$. Those with last $i$ 
	entries increasing in value, $\{\pi \in P^{(0)}_{321}(n):\ 
	\pi_{1}=k,$ last $i \uparrow \}$, are equinumerous with those in 
	which $n$ does not occur among the last $i-1$ entries, $\{\pi \in P^{(0)}_{321}(n):\ 
	\pi_{1}=k,\ n=\pi_{j}$  with  $j \le n-(i-1) \}$.
\end{lemma}
\begin{proof}
	We define a bijection $\phi$ from the first set to the second: $\phi$ 
	is the identity if $n=\pi_{j}$ with $j\le n-i$, otherwise, the 
	hypothesis ``last $i \uparrow$'' implies $n=\pi_{n}$ and $\phi$ 
	interchanges $\pi_{n}$ and $\pi_{n-i+1}$.
\end{proof}
\begin{lemma}
	\label{09}
	\[
	\begin{array}{l}
	\#\{\pi \in P^{(0)}_{321}(n):\ 
	\pi_{1}=k,\ n \textrm{ does not occur among the last $i-1$ entries} 
	\}= \\
	\#\{\t{Dyck $n$-paths:\ first ascent $\ge k$ and last descent $\ge 
	i$} \}.
	\end{array}
	\]
\end{lemma}
\begin{proof}
    Krattenthaler's bijection sends	the first set to the second. 
	This is because $\pi \in P^{(0)}_{321}(n)$ does not have $n$ among 
	the last $i-1$ entries 
	$\Leftrightarrow $ the last location of a record high in $\pi$  is $\le 
	n-(i-1)
	\Leftrightarrow$ the last entry of $\De(record\ high\ 
	positions,\,n+1)$ is $\ge i \Leftrightarrow$ the last descent of the corresponding 
	Dyck path is $\ge i$.
\end{proof}
In particular, using Lemmas \ref{05}\,(i) and \ref{03}\,(ii),
$\#\{\pi \in P^{(0)}_{321}(n):\ $first entry $\ge 2$ and last    
 entry $\le n-1\}=C^{(1)}_{n-2}+C^{(3)}_{n-3}+C^{(5)}_{n-4}=C^{(2)}_{n-2}+C^{(5)}_{n-4}.$

\vspace*{10mm}
{\Large \textbf{4 \quad One 132 Pattern}  } 

As a warmup, we treat this case in full detail. 
\begin{lemma}
	\label{10}
	Suppose $\pi \in P_{132}^{(1)}(n)$ and $a\,c\,b$ is the   unique $132$ 
	pattern in $\pi$. Then $a,c$ are consecutive in  position in 
	$\pi$ and $a,b$ are consecutive in value $($i.e. $b=a+1)$.
\end{lemma}
\begin{proof}
	Suppose some letter, $i$ say, appears between $a$ and $c$ in $\pi$ so 
	that $\pi=\cdots a \cdots i \cdots c \cdots b \cdots $. If $i>b$, 
	then $a\,i\,b$ is another 132 pattern in $\pi$ and if $i<b$, then 
	$i\,c\,b$ is another 132 pattern. Thus there are no possibilities 
	for $i$ and $a,\,c$ are consecutive in position. Now suppose $i$ 
	lies between $a$ and $b$ in value: $a<i<b$. If $i$ occurs before $c$ 
	in $\pi$, then $i\,c\,b$ is a 132, while if $i$ occurs after $c$, 
	then $a\,c\,i$ is a 132. Again there are no possibilities for $i$ 
	and $b=a+1$.
\end{proof}
\begin{lemma}
	\label{11}
	$\#\{\pi\in P_{132}^{(1)}(n): $ the $132$ pattern occupies consecutive 
	positions in $\pi\} = (n-2)C_{n-2}=\binom{2n-4}{n-3}$. 
\end{lemma}
\begin{proof}
	Given $\pi\in P_{132}^{(1)}(n)$ with $a,c,b$ in positions 
	$k,\,k+1,\,k+2$, 
	write $\pi$ as $W_{1}\,a\,c\,b\,W_{2}$. All letters in $W_{1}$ are $>b$ (if 
	$x<b,\ x \in W_{1}$, then $x\,c\,b$ is a 132) and no letter in 
	$W_{2}$ lies between $b$ and $c$ in value (if $b<x<c,\ x \in W_{2}$, 
	then $a\,c\,x$ is a 132). 
	
	We now claim the map $\pi \to (\rho,k)$ with 
    $\rho=\textrm{reduce}(W_{1}\,c\,W_{2})$
	is a bijection to $P^{(0)}_{132}(n-2) \times [n-2]$; since $\v 
	P^{(0)}_{132}(n)\v = C_{n}$, the lemma will follow. To form $\rho,\ 
	a$ and $b$ will be deleted from $\pi$ and 2 subtracted from all 
	letters $>b$. So $\rho$ will have the form $W_{1}'\,c'\,W_{2}'$ with 
	$W_{1}'=W_{1}-2$ (entrywise), $ c'=c-2$ in position $k$, and the 
	number of  successors of $c'$ in $\rho$ that are $<c'$ will be 
	precisely $a$. This permits $a$ and hence $b=a+1$ to be recovered from 
	$(\rho,k)$ and so $\pi$ can be reconstructed by adding 2 to each 
	letter of $\rho$ that is $>a$ and inserting $a,b$ so they occupy 
	positions $k$ and $k+2$.
\end{proof}
\begin{lemma}
	\label{12}
	$\#\{\pi\in P_{132}^{(1)}(n):\ $the $132$ pattern occupies the first, 
	second and last 
	positions in $\pi\}=C_{n-3}$.
\end{lemma}
\begin{proof}
	Here $\pi$ has the form $a\,c\,W_{2}\,b$ and,  
	necessarily,  $(a,c,b)=(n-2,n,n-1)$ (else more than one 132 pattern). 
	Hence the map $\pi \to W_{2}$ is a bijection 
	to the set
 	$P_{132}^{(0)}(n-3)$, counted by $C_{n-3}$.
\end{proof}

Now we count $P_{132}^{(1)}(n)$ by the number $k$ of letters 
between  the $c$ and $b$ of the 132 pattern.
\begin{prop}
	\label{14}
	For $k\in [0,n-3]$, there is a bijection $\phi:\{\pi \in 
	P_{132}^{(1)}(n):\ $the $c$ and $b$ of the $132$ 
	pattern are separated by $k$ other letters$\} \to \{ \rho \in 
	P_{132}^{(1)}(n-k):\ $ the $132$ pattern of $\rho$ occupies consecutive 
	positions in  $\rho\}\,\times\,\{\s \in P_{132}^{(1)}(k+3):\ $ the $132$ pattern 
	of $\s$ occupies the first, 
	second and last positions in $\s\}$.
\end{prop}
\begin{proof}
	Here $ \pi$ has the form 
	$W_{1}\,a\,c\,W_{2}\,b\,W_{3}$ with $k$ letters in the subword 
	$W_{2}$. In fact, $k<a$ and the letters in $W_{2}$ comprise 
	$\{a-1,a-2,\ldots,a-k\}$, else too many 132's. Also, and for the 
	same reason, 
	the letters $1,2,\ldots,a-k-1$ all occur in $W_{3}$.
	These facts imply the map $\phi:\ \pi \to (\rho,\s)$ with 
	$\rho=$ reduce$(W_{1}\,a\,c\,b\,W_{3}) \in P_{132}^{(1)}(n-k)$ 
	and $\s=$  reduce$(a\,c\,W_{2},b) \in P_{132}^{(1)}(k+3)$ is the 
	desired bijection.
		This is because $\rho$ is obtained from 
    $(W_{1}\,a\,c\,b\,W_{3})$ by subtracting $k$ from all letters 
	$\ge a$ and $\s =k\!+\!1\ k\!+\!3\ W_{2}\!-\!(a\!-\!k\!-\!1)\ k\!+\!2$. Hence 
	the 132 pattern in $\rho$ identifies the values of $a,b,c$
	and then knowing $a$, $W_{1}$ and $W_{3}$ can be recovered from 
	$\rho$ as can $W_{2}$ from $\s$. 
\end{proof}
Hence
\begin{eqnarray*}
	\v P_{132}^{(1)}(n)\v & = & \sum_{k=0}^{n-3}\#\{\rho\in 
	P_{132}^{(1)}(n-k):\ a,\,c,\,b\t{ occupy consecutive positions in } 
	\rho\}\times \\
	  & & \phantom{\sum} \#\{\s \in P_{132}^{(1)}(k+3):\ a,\,c,\,b\t{ occupy the first, 
	second and last positions in }\s\}  \\
	 & \overset{(1)}{=} & \sum_{k=0}^{n-3}\binom{2(n-k)-4}{(n-k)-3} \times 
	 C_{(k+3)-3} \\
	   & = & \sum_{k=0}^{n-3}\binom{2n-2k-4}{n-k-3}C_{k} \\
	   & = & \binom{2n-3}{n-3},
\end{eqnarray*}
(1) by Lemmas \ref{11} and \ref{12}, the last equality by Lemma 
\ref{13}. We have established
\begin{theorem}[B\'{o}na \cite{bona98}]
	\label{015}
	$\v P_{132}^{(1)}(n)\v =\binom{2n-3}{n-3}$.
\end{theorem}

\vspace*{10mm}

{\Large \textbf{5 \quad One 321 Pattern}  } 

We wish to show bijectively that $\vert 
P_{321}^{(1)}(n)\vert=C^{(6)}_{n-3}\ 
\big(=\frac{3}{n}\binom{2n}{n-3}\big)$ for $n\ge 3$.  
\begin{lemma}
	\label{016}
	Suppose the $321$ 
	pattern in  $\pi \in P_{321}^{(1)}(n)$ has middle letter $b$. Then the 
	letters preceding $b$ in $\pi$
    consist of the $b-1$ letters  $\{c\} \cup [b-1] \setminus \{a\}$ 
    $($and so $b$ is a fixed point$)$.
\end{lemma}
\begin{proof} The letters preceding $b$ in $\pi$ include $c$ but no 
other letter $> b$ (else more than one $321$ pattern). Among the 
letters following $b$, $a$ is the only one $<b$ (else again more than one 
$321$ pattern). The result follows. 
\end{proof}
\begin{prop}
	Let $b \in [2,n-1]$. Then there is a bijection from $\{\pi \in 
	P_{321}^{(1)}(n):\ $middle letter of $321$ pattern is $b\} \longrightarrow
	\{\rho \in 
	P_{321}^{(0)}(b):\ $last entry $\le b-1\} \times
	\{\s \in 
	P_{321}^{(0)}(n-b+1):\ $first entry $\ge 2\}$.
	\label{bijone321}
\end{prop}
\begin{proof}
	Write $\pi$ as $W_{1}\,b\,W_{2}$ with $c$ occurring in the subword 
	$W_{1}$ and $a$ in $W_{2}$. Set $\rho=$
	reduce$(W_{1}\,a)$ and $\s=$ reduce$(c\,W_{2})$. 
	Note that
	$\rho$ is simply $(W_{1}'\,a)$ where $W_{1}'$ is $W_{1}$ with c 
	replaced by $b$ and, since $\{c,W_{2}\}$
	consists of all letters $>b$ together with $a$, $\s$ is obtained from
	$c\,W_{2}$ by
	replacing $a$ by 1, and subtracting $b-1$ from every other entry.
	It is clear that $\rho,\s$ have the specified
	properties.
	Also the last letter of $\rho$ is $a$ and $\rho$ contains all information 
	about $W_{1}$ except the value 
	of $c$. On the other hand, $\s$ contains the value of $c$ (via its 
	first letter) and all information 
	about $W_{2}$ except the value 
	of $a$. 
    Thus we can uniquely recover $\pi$ 
	from $\rho$ and $\s$ and the map $\pi \to (\rho,\s)$ is a bijection as 
	claimed.

\end{proof}

The sets in the Cartesian product of the preceding 
Proposition are counted by 
$C^{(3)}_{b-2}$ and $C^{(3)}_{n-b-1}$ respectively (Lemma \ref{07}). Hence, counting by the 
middle letter of the 321 pattern, we find
\[
\vert P_{321}^{(1)}(n) \vert = 
\sum_{b=2}^{n-1}C^{(3)}_{b-2}C^{(3)}_{n-b-1}=
\sum_{b=0}^{n-3}C^{(3)}_{b}C^{(3)}_{n-3-b}=C^{(6)}_{n-3}
\]
and we obtain
\begin{theorem}[Noonan \cite{noonan}]  
	\label{18}
    $\vert P_{321}^{(1)}(n) \vert = C^{(6)}_{n-3}$.  
\end{theorem}  
It is observed in \cite{noonan} that  the number of permutations in 
$P_{321}^{(0)}(n+2)$ in which $1,2,3,4,5$ occur in increasing order is 
also $C^{(6)}_{n-3}$. This can be explained as follows: incrementing by 1 all entries 
$\ge 6$ in such a permutation and inserting the letter 6 at the start 
gives a bijection to $\{\pi \in 
P^{(0)}_{321}(n+3):\ \pi_{1}=6\}$, counted by  $C^{(6)}_{n-3}$ (Lemma 
\ref{07}).                                                                                                      

\vspace*{10mm}

{\Large \textbf{6 \quad Two 321 Patterns}  }

For $ P_{321}^{(2)}(n)$ we distinguish two cases according as the two 321 patterns 
have (I)  a common (i.e. the same) middle letter $b$ or (II) distinct 
middle letters. Case I is counted by $b$ itself, case II by the number 
of intervening letters.

\textbf{Case I: Same middle letter}.\quad Here, the two 321 patterns also have a common first 
	letter or 
	a common last letter (or there would be at least four 321 patterns) 
	and by symmetry there are the same number of each type. We count the 
	second type, with patterns $c_{1}ba$ and $c_{2}ba$ so that $\pi$ can 
	be written $W_{1}bW_{2}$ with $c_{1},c_{2}$ occurring in $W_{1}$ 
	and $a$ in $W_{2}$.
\begin{lemma}
	\label{19}
	For $\pi$ as just described, the subword $W_{1}a$ consists of the letters 
	$[b-1]\cup\{c_{1},c_{2}\}$ and $c_{1}c_{2}W_{2}$ consists of 
	$\{a\}\cup [b+1,n]$. Consequently, $\v W_{1} \v =b$ and $b$ is in 
	position $b+1 \in [3,n-1]$.
\end{lemma}	
\begin{proof}
	Similar to that of Lemma 15: anything else would give too many 321 patterns.
\end{proof}	
\begin{prop}
	\label{20}
	Let $b \in [2,n-2]$. Then there is a bijection from $\{\pi \in 
	P^{(2)}_{321}(n):\ $the two $321$ patterns in $\pi$ have $b$ as common 
	middle letter and also a common last letter$\}$ to $\{\rho \in P^{(0)}_{321}(b+1):\ $last 
	entry $\le$ largest entry $(=b\!+\!1) -2\} \times \{\s \in P^{(0)}_{321}
    (n-b+1):\ $ the letter $1$ occurs at or after position $3\}$.
\end{prop}
\begin{proof}
	Write $\pi = W_{1}\,b\,W_{2}$ with $c_{1},c_{2}$ in $W_{1}$ and the 
	common $a$ in $W_{2}$. Set $\rho=$ 
	reduce$(W_{1}\,a)$
    and $\s=$ reduce$(c_{1}\,c_{2}\,W_{2})$. The proof that $\rho,\s$ 
    range over the specified sets and contain enough information to 
    recover $\pi$ is much like that of Proposition \ref{bijone321}.
\end{proof}
\begin{cor}
	\label{21}
    $\#\{\pi \in P^{(2)}_{321}(n):\ $ the two $321$ patterns in $\pi$ 
    have a common middle letter$\}=2C^{(8)}_{n-4}.$
\end{cor}
\begin{proof}
	Assume a common last letter (as noted above, the required count is double this count). 
	The sets in the Cartesian product have sizes $C^{(4)}_{b-2}$ and 
	$C^{(4)}_{n-b-2}$ respectively (Lemma 8). Summing over $b \in [2,n-2]$ gives 
	the result.
\end{proof}

\textbf{Case II: Distinct middle letters}.\quad 
Let $c_{1}\,b_{1}\,a_{1}$ and $c_{2}\,b_{2}\,a_{2}$ be the two 321 
patterns with $c_{1}\le c_{2}$ and write $\pi$ as 
$W_{1}\,b_{1}\,W_{2}\,b_{2}\,W_{3}$.  Clearly, $c_{1}$ occurs in $W_{1}$, and $c_{2}$ in $W_{2}$ 
unless $c_{1}=c_{2}$. Also, $a_{2}$ occurs in $W_{3}$, and $a_{1}$ in 
$W_{2}$ unless $a_{1}=a_{2}$. This seems to give three cases (four if you 
distinguish the order $c_{2}$ and $a_{1}$ occur in $W_{2}$) but we 
can handle them all together. In order not to exceed the quota 
of 321 patterns, the subword $W_{1}$ consists of the 
letters $[b_{1}-1]\cup\{c_{1}\}\backslash\{a_{1}\}, \ W_{2}$ consists 
of $[b_{1}+1,b_{2}-1]\cup\{a_{1},c_{2}\}\backslash\{c_{1},a_{2}\}$ 
(in the multiset sense, since $a_{1}$ may $=a_{2}$ or $c_{1}$ may 
$=c_{2}$), and $W_{3}$ consists of $[b_{2}+1,n] \cup 
\{a_{2}\}\backslash \{c_{2}\}$. In any case $\v W_{1}\v=b_{1} 
-1,\ \v W_{2}\v= b_{2}-b_{1}-1,\ \v W_{3}\v = n-b_{2}$.
\begin{prop}
	\label{22}
	Let $k \in [0,n-4]$. Then there is a bijection from $\{\pi \in 
	P^{(2)}_{321}(n):\ $the two $321$ patterns in $\pi$ have distinct 
	middle letters separated by $k$ other letters$\}$ to $\{\rho \in 
	P^{(1)}_{321}(n-k-1)\} \times \{\s \in P^{(0)}_{321}
    (k+2):\ $ first entry of $\s\ge $ smallest entry $+1$, last entry 
    of $\s\le $ largest entry $-1\}$.
\end{prop}
\begin{proof}
	Write $\pi = W_{1}\,b_{1}\,W_{2}\,b_{2}\,W_{3}$.  Set $\rho=$ 
	reduce$(W_{1}'\,b_{1}\,W_{3}')$ where $W_{1}'$ is $W_{1}$ with 
	$c_{1}$ replaced by $c_{2}$ (this includes the case $c_{1}=c_{2}$ in 
	which case $W_{1}'=W_{1}$) and $W_{3}'$ is $W_{3}$ with 
	$a_{2}$ replaced by $a_{1}$. Set $\s=$ 
	reduce$(c_{1}\,W_{2}\,a_{2}).$ Then $\rho,\s$ lie in the specified 
	sets. 
	
	To get $\rho$ from $W_{1}'\,b_{1}\,W_{3}'$, all letters $\le 
	b_{1}$ are left intact and $b_{2}-b_{1}$ is subtracted from all 
	letters $>b_{2}$ (there are no letters in $[b_{1},b_{2}-1]$.) Now 
	$W_{1}$ ( and $W_{1}'$) has only one letter $>b_{1}$ and so the 
	position of $c_{1}$ in $\pi$ is the position of the ``$c$'' in the 
	321 pattern of $\rho$. Adding $b_{2}-b_{1}$ to this ``$c$'' gives 
	the value of $c_{2}$. The ``$a$'' of the pattern is $a_{2}$. The 
	other letters of $\rho$ yield both position and value of all letters 
	in $W_{1}$ and $W_{3}$ other than $c_{1}$ and $a_{2}$. 
	
	To see what $\s$ is, note that $c_{1}\,W_{2}\,a_{2}$ consists of 
	$\{a_{1}\}\cup [b_{1}+1,b_{2}-1]\cup\{c_{2}\}$. Hence reducing will 
	send $a_{1}$ to 1, $c_{2}$ to $k+2$ (since there are $k$ letters in 
	$[b_{1}+1,b_{2}-1])$, and all other letters will be decremented by 
	$b_{1}-1$. So the position of 1 in $\s$ yields the position of 
	$a_{1}$ in $c_{1}\,W_{2}\,a_{2}$ (if 1 is in the last position, 
	then $a_{1}=a_{2}$). Likewise, the position of $k+2$ in $\s$ yields 
	the position of $c_{2}$ in $c_{1}\,W_{2}\,a_{2}$ (if $k+2$ is the 
	first letter, then $c_{2}=c_{1}$). Adding $b_{1}-1$ to the other 
	letters in $\s$ yields positions and values of all letters in 
	$c_{1}\,W_{2}\,a_{2}$ other than $c_{2}$ and $a_{1}$.
	
	Summarizing, $\rho$ gives the nonpattern letters in $W_{1}$ and 
	$W_{3}$, $\s$ those in $W_{2}$; $\rho$ gives the position of $c_{1}$ 
	and $a_{2}$ in $\pi$ and the values of $c_{2}$ and $a_{1}$ while 
	$\s$ gives the positions of $a_{1}$ and $c_{2}$, tells if 
	$a_{1}=a_{2}$ or $c_{2}=c_{1}$ and if not, gives the values of 
	$a_{2}$ and $c_{1}$. So $\rho$ and $\s$ contain just enough information 
	to recover $\pi$ and the map $\pi \to (\rho,\s)$ is a 
	bijection as claimed.
\end{proof}
\begin{cor}
	\label{23}
	$\#\{\pi \in P^{(2)}_{321}(n):\ $ the two $321$ patterns in $\pi$ 
    have distinct middle letters$\}=C^{(8)}_{n-4}+C^{(11)}_{6}.$
\end{cor}
\begin{proof}
	The sizes of the sets in the Cartesian product of the Proposition 
	are $C^{(6)}_{n-k-4}$ (Theorem \ref{18}) and $C^{(2)}_{k}+C^{(5)}_{k-2}$ 
	(Lemma \ref{09}) and summing 
	over $k \in [0,n-4]$ gives the result.
\end{proof}

Combining Cases I and II, we have a result conjectured in \cite{zeilnoonan} 
and proved in \cite{fulmek02}.
\begin{theorem}
	\label{24}
	$\v P^{(2)}_{321}(n)\v=3C^{(8)}_{n-4}+C^{(11)}_{n-6}. $
\end{theorem}

\newpage

{\Large \textbf{7 \quad Three and Four 321 Patterns}  } 

For these situations, counts are also needed of subsets of 
$P_{321}^{(1)}$ and $P_{321}^{(2)}$ satisfying certain conditions on 
the values of the first/last letters and/or the positions of the 
max/min letters. The 
computations are laid out in the following tables (considerable 
details are left to the reader). The expressions obtained are by no 
means unique. For example, the total in the first table below, 
$2C^{(6)}_{n-4}+C^{(9)}_{n-6}$ could also be rendered 
$C^{(5)}_{n-4}+C^{(7)}_{n-4}$ using Lemma 3. Denote the 321 
patterns $c_{1}b_{1}a_{1},c_{1}b_{1}a_{1},\ldots$ ordered primarily 
by increasing $b$, secondarily by increasing $c$ and finally by 
increasing $a$. In the tables, $\Delta C^{(k)}_{n}$ means 
$C^{(k)}_{n}-C^{(k)}_{n-1}$.
The counts in Cases 2 and 3 for three 321 patterns need to be doubled 
due to symmetry considerations. Likewise for Cases 2,\,4 and 5 for 
four 321 patterns.

\begin{center}
\begin{tabular}{|c|c|}  \hline
$\{\pi\in P_{321}^{(1)}(n):\ $ last 2 $\uparrow$\} \\ \hline
\end{tabular}\\
\vspace*{2mm}
\begin{tabular}{|c|c|}  \hline
 Case & there is only one case\\ \hline \hline
pattern type   &$\pi=W_{1}\,b\,W_{2}$          \\ \hline
position of $b$& $b$                                  \\ \hline  \hline
$\rho=$       & reduce$(W_{1}a)$                 \\ \hline
ranges over   & $P_{321}^{(0)}(b)$                       \\  \hline       
subject to    & last $\ne$ max                           \\  \hline            
a set of size &$C^{(3)}_{b-2}$             \\ \hline \hline
        
$\sigma=$      & reduce$(cW_{2})$       \\  \hline
ranges over   & $P_{321}^{(0)}(n-b+1)$             \\  \hline       
subject to    & first $\ge 2$, last 2 $\uparrow$   \\  \hline
a set of size & $2C^{(3)}_{n-b-2}+C^{(6)}_{n-b-4}$         \\ \hline\hline
sum over      & $b\in[2,n-2]$                       \\ \hline
  total    & $2C^{(6)}_{n-4}+C^{(9)}_{n-6}$ \\ \hline
\end{tabular}

\end{center}

\newpage

\begin{center}
\begin{tabular}{|c|}  \hline
$\{\pi\in P_{321}^{(2)}(n):\ $ last 2 $\uparrow$\} \\ \hline
\end{tabular}
\end{center}
\vspace*{-5mm}
\begin{tabular}{|c|c|c|c|}  \hline
 Case & 1 & 2 & 3 \\ \hline \hline
pattern type   & distinct $b$'s                                  & common $b$'s     & common $b$'s \\ 
               &$\pi=W_{1}\,b_{1}\,\un{W_{2}}_{k}\,b_{2}W_{3}$   & common  $a$'s    & common  $c$'s  \\ \hline
position of $b$& both are fixed points                           & $b+1$            & $b-1$ \\ \hline  \hline
$\rho=$       & reduce$(W_{1}'\,b_{1}\,W_{3}')$            & reduce$(W_{1}a)$       & reduce$(W_{1}a_{1}a_{2})$  \\ 
              & $W_{1}'=W_{1}/.\:c_{1}\leftarrow c_{2}$    &                        &    \\ 
              & $W_{3}'=W_{3}/.\:a_{2}\leftarrow a_{1}$    &                        &    \\ \hline
ranges over   & $P_{321}^{(1)}(n-k-1)$                     & $P_{321}^{(0)}(b+1)$   & $P_{321}^{(0)}(b)$     \\  \hline       
subject to    & last 2 $\uparrow$                          & last $\le$ max$-2$     & max in pos $\le -3$    \\  \hline
             
a set of size &$2C^{(6)}_{n-k-5}+C^{(9)}_{n-k-7}$          & $C^{(4)}_{b-2}$        & $C^{(4)}_{b-3}$         \\ \hline \hline
        
$\sigma=$      & reduce$(c_{1}W_{2}a_{2})$      &  reduce$(c_{1}c_{2}W_{2})$        & reduce$(c\,W_{2})$      \\  \hline
ranges over   & $P_{321}^{(0)}(k+2)$            & $P_{321}^{(0)}(n-b+1)$            & $P_{321}^{(0)}(n-b+2)$  \\  \hline       
subject to    & first $>$ min, last $<$ max     & 1 weakly after pos 3              & first $\ge 3$   \\            
              &                                 & and last 2 $\uparrow$             & and last 2 $\uparrow$    \\  \hline
a set of size & $C^{(2)}_{k}+C^{(5)}_{k-2}$     & $C^{(6)}_{n-b-4}+C^{(4)}_{n-b-3}+$   & $C^{(6)}_{n-b-3}+C^{(4)}_{n-b-2}$  \\ 
              &                                 & $C^{(2)}_{n-b-2}-C^{(0)}_{n-b-2}$    & $+C^{(2)}_{n-b-1}$     \\ \hline\hline
sum over      & $k\in[0,n-5]$                   & $b\in[2,n-3]$                        & $b\in[3,n-1]$           \\ \hline
case total    & $2C^{(8)}_{n-5}+3C^{(11)}_{n-7}+C^{(14)}_{n-9}$& $C^{(10)}_{n-6}+C^{(8)}_{n-5}$ 
                                                    & $C^{(10)}_{n-6}+C^{(8)}_{n-5}+C^{(6)}_{n-4}$  \\  
              &                & $+C^{(6)}_{n-4}-C^{(4)}_{n-4}$        &        \\ \hline
\end{tabular}
\vspace*{-5mm}
\begin{center}
\begin{tabular}{|c|}  \hline
Grand Total:\ $4C^{(8)}_{n-5}+3C^{(11)}_{n-7}+C^{(14)}_{n-9}+2C^{(10)}_{n-6}+2C^{(6)}_{n-4}-C^{(4)}_{n-4}$ \\ \hline
\end{tabular}
\end{center}

\newpage

\begin{center}
\begin{tabular}{|c|}  \hline
Three 321 patterns \\ \hline
\end{tabular}\\
\vspace*{5mm}
$\pi=W_{1}\,b\,W_{2}$ \\
\end{center}
\vspace*{-5mm}
\begin{tabular}{|c|c|c|c|}  \hline
 Case & 1 & 2 & 3 \\ \hline \hline
pattern type   & $b_{1}\le b_{2}<b_{3}\equiv b$ & $b_{1}< b_{2}=b_{3}
                                                         \equiv b,\ c_{2}=c_{3}\equiv c$ & common $b$'s and $a$'s  \\ \hline
position of $b$& $b$                      & $b-1$                                        & $b+2$ \\ \hline  \hline
$\rho=$        & reduce$(W_{1}a_{3})$     & reduce$(W_{1}a_{2}a_{3})$                    & reduce$(W_{1}a)$  \\ \hline
ranges over   & $P_{321}^{(2)}(b)$        & $P_{321}^{(1)}(b)$                           & $P_{321}^{(0)}(b+2)$ \\  \hline       
subject to    & last $\ne$ max            & last two $\uparrow$ and max occurs earlier   & last $\le $ max$-3$  \\  \hline       
a set of size & $\Delta P_{321}^{(2)}(b)$ & $C^{(6)}_{b-4}+C^{(9)}_{b-6}$                & $C^{(5)}_{b-2}$  \\ \hline \hline
        
$\sigma=$      & reduce$(c_{3}W_{2})$     &  reduce$(cW_{2})$                            & reduce$(c_{1}c_{2}c_{3}W_{2})$   \\  \hline
ranges over   & $P_{321}^{(0)}(n-b+1)$    & $P_{321}^{(0)}(n-b+2)$                       & $P_{321}^{(0)}(n-b+1)$ \\  \hline       
subject to    & $\si_{1}\ge 2$            &  $\si_{1}\ge 3$                              & 1 weakly after pos 4   \\  \hline          
a set of size & $C^{(3)}_{n-b-1}$         & $C^{(4)}_{n-b-1}$                            & $C^{(5)}_{n-b-3}$  \\ \hline \hline
sum over      & $b\in[4,n-1]$             & $b\in[4,n-1]$                                & $b\in[2,n-3]$  \\ \hline
case total    & $\Delta(3 C^{(11)}_{n-5}+C^{(14)}_{n-7})$ & $C^{(10)}_{n-5}+C^{(13)}_{n-7}$ & $ C^{(10)}_{n-5}$  \\ \hline 

\end{tabular}
\vspace*{-5mm}
\begin{center}
\begin{tabular}{|c|}  \hline
Grand Total:\ $7C^{(10)}_{n-5}+6C^{(13)}_{n-7}+C^{(16)}_{n-9}$ \\ \hline
\end{tabular}
\end{center}

\newpage
 
\vspace*{-20mm}
\begin{center}
\begin{tabular}{|c|}  \hline
Four 321 patterns \\ \hline
\end{tabular}

\begin{tabular}{|c|c|c|c|}  \hline
 Case          & 1 & 2 & 3 \\ \hline \hline
pattern type   & $b_{3}<b_{4}\equiv b$       & $b_{2}< b_{3}=b_{4}\equiv b,\ c_{3}=c_{4}\equiv c$ & 4321 pattern present \\ \hline											                                			  					  
position of $b$& $b$                        & $b-1$ & $c,b$ adjacent, $c=b+1$\\ 
                &                           &       & $\pi=W_{1}\,c\,b\,W_{2}$    \\ \hline  \hline  
$\rho=$        & reduce$(W_{1}a_{4})$       & reduce$(W_{1}a_{3}a_{4})$  & $W_{1}cW_{2}$   \\ \hline
ranges over    & $P_{321}^{(3)}(b)$         & $P_{321}^{(2)}(b)$     &  $P_{321}^{(1)}(n-1)$ \\  \hline       
subject to     & last $\ne$ max             & last two $\uparrow$ and max occurs earlier & ---  \\  \hline       
a set of size  & $\Delta P_{321}^{(3)}(b)$  & 
$C^{(8)}_{b-5}+2C^{(11)}_{b-7}+C^{(14)}_{b-9}+$  &$C^{(6)}_{b-4}$ \\ 
               &                            & 
               $2C^{(10)}_{b-6}+2C^{(6)}_{b-4}-C^{(4)}_{b-4}$ &  \\ \hline   \hline       
$\sigma=$      & reduce$(c_{4}W_{2})$       & reduce$(cW_{2})$ & not needed  \\  \hline
ranges over    & $P_{321}^{(0)}(n-b+1)$     & $P_{321}^{(0)}(n-b+2)$  & \\  \hline       
subject to     & $\si_{1}\ge 2$             & $\si_{1}\ge 3$ &    \\  \hline          
a set of size  & $C^{(3)}_{n-b-1}$          & $C^{(4)}_{n-b-1}$ &  \\ \hline \hline
sum over       & $b\in[5,n-1]$              & $b\in[4,n-1]$  &  \\ \hline

case total     & $\Delta(7 C^{(13)}_{n-6}+6C^{(16)}_{n-8}$ 
                                            & $C^{(12)}_{n-6}+2C^{(15)}_{n-8}+C^{(18)}_{n-10}$& $C^{(6)}_{n-4}$  \\  
               & $+C^{(19)}_{n-10})$        & $+2C^{(14)}_{n-7}+2C^{(10)}_{n-5}-C^{(8)}_{n-5}$& \\ \hline
\end{tabular}

\begin{tabular}{|c|c|c|c|}  \hline
 Case & 4 & 5 & 6 \\ \hline \hline
pattern type   & $b_{1}< b_{2}=b_{3}=b_{4}\equiv b$ & common $b$'s                   & common $b$'s \\ 
               & $c_{2}=c_{3}=c_{4}\equiv c$        & common $a$'s                   & $c_{1}=c_{2}<c_{3}=c_{4}$ \\ \hline
position of $b$& $b-2$                              & $b+3$                          & $b$ \\ \hline  \hline
$\rho=$       & reduce$(W_{1}a_{2}a_{3}a_{4})$      & reduce$(W_{1}a)$               & reduce$(W_{1}a_{2}a_{3})$  \\ \hline
ranges over   & $P_{321}^{(1)}(b)$                  & $P_{321}^{(0)}(b+3)$           & $P_{321}^{(0)}(b+1)$ \\  \hline       
subject to    & last 3 $\uparrow$ and               & last $\le$ max$-4$             & two largest   \\  
              &  max occurs earlier                 &                                & 3 or more from end \\  \hline  
a set of size &$C^{(7)}_{b-5}+C^{(10)}_{b-7}$       & $C^{(6)}_{b-2}$                & $C^{(5)}_{b-3}$  \\ \hline \hline
        
$\sigma=$      & reduce$(cW_{2})$              &  reduce$(c_{1}c_{2}c_{3}c_{4}W_{2})$     & reduce$(c_{3}c_{4}W_{2})$   \\  \hline
ranges over   & $P_{321}^{(0)}(n-b+3)$         & $P_{321}^{(0)}(n-b+1)$                   & $P_{321}^{(0)}(n-b+2)$ \\  \hline       
subject to    & $\si_{1}\ge 4$                 & 1 weakly after pos 5                     & 1 and 2 weakly after pos 3   \\  \hline          
a set of size & $C^{(5)}_{n-b-1}$              & $C^{(6)}_{n-b-4}$                        & $C^{(5)}_{n-b-2}$  \\ \hline \hline
sum over      & $b\in[5,n-1]$                  & $b\in[2,n-4]$                            & $b\in[3,n-2]$  \\ \hline
case total    & $C^{(12)}_{n-6}+C^{(15)}_{n-8}$& $C^{(12)}_{n-6}$                         & $C^{(10)}_{n-5}$        \\ \hline 

\end{tabular}
\end{center}
\begin{tabular}{|c|}  \hline
Grand Total:\ 
$13C^{(12)}_{n-6}+19C^{(15)}_{n-8}+9C^{(18)}_{n-10}+C^{(21)}_{n-12}+4C^{(14)}_{n-7}
+5C^{(10)}_{n-5}+C^{(6)}_{n-4}-2C^{(8)}_{n-5}$ \\ \hline
\end{tabular}

\newpage

{\Large \textbf{8 \quad 132 Patterns in 123-Avoiding Permutations}  } 

Since we are henceforth dealing with 123-avoiding permutations, let $P^{(i)}(n)$ denote the set 
of 123-avoiding permutations on $[n]$ that contain exactly $i\ 132$ 
patterns. Here we compute $\v P^{(i)}(n) \v$ for $1\le i \le 4$. We 
distinguish cases according to the rightmost $c$ among the $i$ 132 patterns. 
Say this $c$ does single duty if it appears as the ``$c$'' of only one 
132 pattern, double duty if it is the ``$c$'' of two 
132 patterns, and so on. To start the ball rolling, we have the 
following count for $P^{(0)}(n)$. A proof is 
included for completeness.
\begin{theorem}[\cite{simion-schmidt}]	The number of 
     permutations on $[n]$ that avoid $123$ and $132$  is $2^{n-1}$.
\end{theorem}
\begin{proof}
	For $\pi$ on $[n]$ to avoid 123 and 132 patterns, $\pi_{1}$ must be 
	$n$ or $n-1$ for if it were $\le n-2$ it would have two larger 
	successors, forcing $\pi_{1}$ to initiate a proscribed pattern. 
	Similarly, $\pi_{k}\ge n-k$ for all $1\le k \le n-1$ and, conversely, this condition 
	guarantees that $\pi$ avoids 123 and 132. The number of such permutations 
	is the permanent of the lower Hessenberg matrix
	\[
\bordermatrix{
  & n  & n-1   & n-2    &  \ldots  & \ldots   &  1      \cr
1     & 1  & 1  &   0  &   0  &   \ldots   &  0   \cr
2     & 1  & 1   & 1   &   0  & \ldots  &    0  \cr
3     & 1  & 1   & 1    &  1   &      &    \vdots  \cr
\vdots& \vdots  & \vdots   &  \vdots   &  \ddots   &  \ddots    &   0   \cr
\vdots& 1  & 1   &  1   &   \ldots   &   1  &    1  \cr
n     & 1  & 1  &1    &  \ldots    &   1  & 1     \cr
}
\]
which is well known to be $2^{n-1}$.	
\end{proof}
We now proceed to count $P^{(i)}(n)$ for $1\le i \le 4$, considering 
separate cases according to pattern overlap and summarizing 
computations in tables similar to those in the last section. The 
reader is left to verify details.

\vspace*{10mm}

{\large \textbf{No 123, One 132 Pattern}  } 

By Lemma \ref{10}, $\pi$ has the form $W_{1}\,a\,c\,W_{2}$ with $b \in 
W_{2}$ and $b=a+1$. Also, $W_{1}$ comprises the letters 
$[b+1,n]\backslash \{c\}$ and $W_{1}$ comprises $[1,b]\backslash \{a\}$. 
For $b \in [2,n-1]$, the map $\pi \to (\rho,\s)$ is a bijection 
where $\rho=$ reduce($W_{1}c)$ ranges over $P^{(0)}(n-b)$, a set of 
size $2^{n-b-1}$, and $\s=$ reduce$(W_{2})$  ranges over $P^{(0)}(b-1)$, a set of 
size $2^{b-2}$. Summing over $b \in [2,n-1]$, we have 
\begin{theorem}[Robertson \cite{aaron2}]
	$\v P^{(1)}(n)\v=\binom{n-2}{1}2^{n-3}$.
\end{theorem}

\vspace*{10mm}

{\large \textbf{No 123, Two 132 Patterns}  } 

Here the two cases are (1) $c_{2}$ does single duty, that is, $c_{1}\ne c_{2}$, (2) $c_{2}$ does double duty, that is,
$c_{1}=c_{2}$.  In 
case (1) $\pi$ has the form $W_{1}\,a_{2}\,c_{2}\,W_{2}$ with $b_{1}\in W_{1} \cup \{c_{2}\},\ a_{1},c_{1} \in 
W_{1},\ b_{2} \in W_{2}$ and with $ a_{2},b_{2}$ consecutive in value and $b_{2}<b_{1}$ (else $a_{1}b_{1}b_{2}$ is a 
123). Case (2) 
forces common $b$'s as well as common $c$'s, and $\pi$ has the form  
$W_{1}\,a_{1}\,a_{2}\,c\,W_{2}$ with $b= b_{1}=b_{2} \in W_{2}$ and 
$a_{2},a_{1},b$ consecutive in value.

\begin{center}
\begin{tabular}{|c|}  \hline
$\{\pi\in P^{(2)}(n)\}$ \\ \hline
\end{tabular}
\\ 
\vspace*{2mm}
\begin{tabular}{|c|c|c|}  \hline
 Case & 1 & 2  \\ \hline \hline
pattern type          & $c_{1}\ne c_{2}$                  & $c_{1}= c_{2}$ \\  \hline
$\pi$ has form        &$W_{1}\,a_{2}\,c_{2}\,W_{2}$       & $W_{1}\,a_{1}\,a_{2}\,c\,W_{2}$  \\ \hline
  $a,c,b=$            & $a_{2},c_{2},b_{2}$               & obvious \\ \hline  \hline
$W_{1}$ comprises     & $[b+1,n]\b \{c\}$                 & $[b+1,n]\b \{c\}$        \\	\hline
$W_{2}$ comprises     & $[1,b]\b \{a\} $                  & $[1,b]\b \{a_{1},a_{2}\} $       \\ \hline	\hline
$\rho=$               & reduce$(W_{1}\,c)$                & reduce$(W_{1}\,c)$  \\ \hline
ranges over           & $P^{(1)}(n-b)$                    & $P^{(0)}(n-b)$     \\  \hline       
a set of size         &$\binom{n-b-2}{1}2^{n-b-3}$        & $2^{n-b-1}$         \\ \hline \hline        
$\sigma=$             & reduce$(W_{2})$                   & reduce$(W_{2})$      \\  \hline
ranges over           & $P^{(0)}(b-1)$                    & $P^{(0)}(b-2)$  \\  \hline       
a set of size         & $2^{b-2}$                         & $2^{b-3}$  \\  \hline\hline
sum over              & $b\in[2,n-3]$                     & $b\in[3,n-1]$           \\ \hline
case total            & $\binom{n-3 }{2}2^{n-5}$          & $\binom{n-3}{1}2^{n-4}$    \\ \hline                                   
\end{tabular}
\\ 
\vspace*{2mm} 
\begin{tabular}{|c|}  \hline
Grand Total:\ $\binom{n-3}{1}2^{n-4}+\binom{n-3}{2}2^{n-5}$ \\ \hline
\end{tabular}
\end{center}

\newpage 

{\large \textbf{No 123, Three 132 Patterns}  } 

Here there are four cases: (1) $c_{3}$ does single duty and $b_{3}$ is the only pattern letter following $c_{3}$,
(2) $c_{3}$ does single duty and $b_{2}=b_{3}$, (3) $c=c_{2}=c_{3}$ does double duty, 
(4) $c=c_{1}=c_{2}=c_{3}$ does triple duty. 

In case (1) $\pi = W_{1}\,a\,c\,W_{2}$ with $a=a_{3},\ c=c_{3},\ b=b_{3} \in W_{2}$
and $a,b$ consecutive in value. 

In case (2) only four distinct letters 
are involved among the three 132 patterns (their reduced form is 
1432 with patterns $1\,4\,3, 1\,4\,2, 1\,3\,2$) and $\pi = W_{1}\,a\,c\,b\,W_{2}$ with $a=a_{1}=a_{2}=a_{3},\ 
c=c_{1}=c_{2},\ b=b_{1}=c_{3}$ and $b_{2}=b_{3}\in W_{2}$. Also
$a,b_{2},b$ consecutive in value.

In case (3)  $\pi = W_{1}\,a_{2}\,a_{3}\,c\,W_{2}$ with $
b=b_{2}=b_{3} \in W_{2}$
and $a_{3},a_{2},b$ consecutive in value.

In case (4)  $\pi = W_{1}\,a_{1}\,a_{2}\,a_{3}\,c\,W_{2}$ with $
b=b_{1}=b_{2}=b_{3} \in W_{2}$
and all four of $a_{3},a_{2},a_{1},b$ consecutive in value.

\begin{center}
\begin{tabular}{|c|}  \hline
$\{\pi\in P^{(3)}(n)\}$ \\ \hline
\end{tabular}
\\ 
\vspace*{2mm}
\begin{tabular}{|c|c|c|}  \hline
 Case & 1 & 2  \\ \hline \hline
pattern type          & $c_{3}$ single duty, $b_{2}\ne b_{3}$   &$c_{3}$ single duty, $b_{2}= b_{3}$\\   \hline
$\pi$ has form        &$W_{1}\,a\,c\,W_{2}$               & $W_{1}\,a\,c\,b\,W_{2}$  \\ \hline
  $a,c,b=$            & $a_{3},c_{3},b_{3}$               & $a_{1},c_{1},b_{1}$ \\ \hline  \hline
$W_{1}$ comprises     & $[b+1,n]\b \{c\}$                 & $[b+1,n]\b \{c\}$        \\	\hline
$W_{2}$ comprises     & $[1,b]\b \{a\} $                  & $[1,b-1]\b \{a=b-2\} $       \\ \hline	\hline
$\rho=$               & reduce$(W_{1}\,c)$                & reduce$(W_{1}\,c)$  \\ \hline
ranges over           & $P^{(2)}(n-b)$                    & $P^{(0)}(n-b)$     \\  \hline       
a set of size         & $\binom{n-b-3}{1}2^{n-b-4}+
                         \binom{n-b-3}{2}2^{n-b-5}$       & $2^{n-b-1}$         \\ \hline \hline        
$\sigma=$             & reduce$(W_{2})$                   & reduce$(W_{2})$      \\  \hline
ranges over           & $P^{(0)}(b-1)$                    & $P^{(0)}(b-2)$  \\  \hline       
a set of size         & $2^{b-2}$                         & $2^{b-3}$  \\  \hline\hline
sum over              & $b\in[2,n-4]$                     & $b\in[3,n-1]$           \\ \hline
case total            & $\binom{n-4}{2}2^{n-6}+
                         \binom{n-4}{3}2^{n-7}$           &  $\binom{n-3}{1}2^{n-4}$    \\ \hline                                   
\end{tabular}
\\ 
\newpage 
\begin{tabular}{|c|}  \hline
$\{\pi\in P^{(3)}(n)\}$ continued \\ \hline
\end{tabular}
\\ 
\vspace*{2mm}
\begin{tabular}{|c|c|c|}  \hline
 Case & 3 & 4 \\ \hline \hline
pattern type          & $c_{3}$ double duty               & $c_{3}$ triple duty  \\  \hline
$\pi$ has form        &$W_{1}\,a_{2}\,a_{3}\,c\,W_{2}$    & $W_{1}\,a_{1}\,a_{2}\,a_{3}\,c\,W_{2}$  \\ \hline
  $a,c,b=$            & $a_{3},c_{3},b_{3}$               & obvious \\ \hline  \hline
$W_{1}$ comprises     & $[b+1,n]\b \{c\}$                 & $[b+1,n]\b \{c\}$        \\	\hline
$W_{2}$ comprises     & $[1,b]\b \{a_{2},a_{3}\} $        & $[1,b]\b \{a_{1},a_{2},a_{3}\} $       \\ \hline	\hline
$\rho=$               & reduce$(W_{1}\,c)$                & reduce$(W_{1}\,c)$  \\ \hline
ranges over           & $P^{(1)}(n-b)$                    & $P^{(0)}(n-b)$     \\  \hline       
a set of size         &$\binom{n-b-2}{1}2^{n-b-3}$        & $2^{n-b-1}$         \\ \hline \hline        
$\sigma=$             & reduce$(W_{2})$                   & reduce$(W_{2})$      \\  \hline
ranges over           & $P^{(0)}(b-2)$                    & $P^{(0)}(b-3)$  \\  \hline       
a set of size         & $2^{b-3}$                         & $2^{b-4}$  \\  \hline\hline
sum over              & $b\in[3,n-3]$                     & $b\in[4,n-1]$           \\ \hline
case total            & $\binom{n-4}{2}2^{n-6}$           & $\binom{n-4}{1}2^{n-5}$    \\ \hline                                   
\end{tabular}
\\ 
\vspace*{2mm} 
\begin{tabular}{|c|}  \hline
Grand Total:\ $\binom{n-3}{1}2^{n-4}+\binom{n-3}{2}2^{n-5}+\binom{n-4}{3}2^{n-7}$ \\ \hline
\end{tabular}
\end{center}

\vspace*{10mm}

{\large \textbf{No 123, Four 132 Patterns}  } 

The 7 cases to be considered here are

Case 1.\quad $c_{4}$ single duty, $b_{4}$ is the only pattern letter occurring 
after $c_{4}$. Here $\pi= W_{1}\,a\,c\,W_{2}$ with 
$a=a_{4},\,c=c_{4},\,b=b_{4}\in W_{2}$ and $b=a+1$.

Case 2.\quad $c_{4}$ single duty, $b_{3}=b_{4}$, only two distinct $c$'s.
Here there are only 5 distinct pattern letters and the pattern overlap 
is as in the figure below with $a,c,b$ consecutive in position (all 
entries in the same column coincide).

\begin{center}
\begin{tabular}{|l c| c | c | c |c l|} \hline 
	$\phantom{a_{1}}$ & $\phantom{a_{1}}$ & $a$ & $c$ & $b$ & $\phantom{a_{1}}$ & 
	$\phantom{a_{1}}$\\ \hline
	 $a_{1}$ & & & $c_{1}$ & $b_{1}$ & & \\
	       & & $a_{2}$ & $c_{2}$ & $b_{2}$ & & \\
		   & & $a_{3}$ & $c_{3}$ & & & $b_{3}$ \\
		   & & $a_{4}$ & & $c_{4}$ & &$b_{4}$ \\ \hline
\end{tabular}
\end{center}
	
Writing $\pi= W_{1}\,a\,c\,b\,W_{2}$ we must also have all four of 
$a,b_{3},a_{1},b$ consecutive in value: $a=b-3,\ b_{3}=b-2,\ 
a_{1}=b-1$.

Case 3.\quad $c_{4}$ single duty, $b_{3}=b_{4}$, three distinct 
$c$'s (four distinct $c$'s is not possible).
Here there are 6 distinct pattern letters and the pattern overlap 
is as in the figure below again with $a,c,b$ consecutive in position.

\begin{center}
\begin{tabular}{|l c c c| c | c | c |c l|} \hline 
 $\phantom{a_{1}}$ &	$\phantom{a_{1}}$ & $\phantom{a_{1}}$ & $\phantom{a_{1}}$ & $a$ & $c$ & $b$ & $\phantom{a_{1}}$ & 
	$\phantom{a_{1}}$\\ \hline
	 $a_{1}$ & & $c_{1}$ &  & & $b_{1}$ &  & & \\
	    & &   & & $a_{2}$ & $c_{2}$ & $b_{2}$ & & \\
		& &   & & $a_{3}$ & $c_{3}$ & & & $b_{3}$ \\
		& &   & & $a_{4}$ & & $c_{4}$ & &$b_{4}$ \\ \hline
\end{tabular}
\end{center}

Writing $\pi= W_{1}\,a\,c\,b\,W_{2}$, we have
$b_{3}\,(=b_{4}),\,b$ consecutive in value: $b_{3}=b-1$.

Case 4.\quad $c_{4}$ double duty, only 5 distinct pattern letters.
The pattern overlap now has the first four pattern letters consecutive 
in position. 
\begin{center}
\begin{tabular}{| c| c | c | c |c l|} \hline 
      $b-1$ & $c_{2}$ & $b-2$   & $c_{4}$ & $\phantom{a_{1}}$ & $b$\\ \hline
	$a_{1}$ & $c_{1}$ &         & $b_{1}$ &                   &          \\
	$a_{2}$ & $c_{2}$ &         &         &                   & $b_{2}$  \\
	$a_{3}$ &         &         & $c_{3}$ &                   & $b_{3}$  \\
            &         & $a_{4}$ & $c_{4}$ &                   & $b_{4}$  \\ \hline
\end{tabular}
\end{center}

Case 5.\quad $c_{4}$ double duty, more than 5 distinct pattern 
letters. Here $c=c_{3}=c_{4},\ b=b_{3}=b_{4}$ and 
$\pi=W_{1}\,a_{3}\,a_{4}\,c\,W_{2}$ with $a_{4},a_{3},b$ consecutive 
in value: $a_{4}=b-2,\,a_{3}=b-1$.

Case 6.\quad $c=c_{2}=c_{3}=c_{4}$ triple duty. Here 
$\pi=W_{1}\,a_{2}\,a_{3}\,a_{4}\,c\,W_{2}$ with 
$b=b_{2}=b_{3}=b_{4}\in W_{2}$ and $a_{4},a_{3},a_{2},b$ consecutive 
in value.

Case 7.\quad $c=c_{1}=c_{2}=c_{3}=c_{4}$ quad duty. Here 
$\pi=W_{1}\,a_{1}\,a_{2}\,a_{3}\,a_{4}\,c\,W_{2}$ with 
$b=b_{1}=b_{2}=b_{3}=b_{4}\in W_{2}$ and $a_{4},a_{3},a_{2},a_{1},b$ consecutive 
in value.
 
The computation tables follow.

\newpage 
\begin{center}
\begin{tabular}{|c|}  \hline
$\{\pi\in P^{(4)}(n)\}$ \\ \hline
\end{tabular}
\\ 
\vspace*{2mm}
\begin{tabular}{|c|c|c|c|}  \hline
 Case & 1 & 2 & 3 \\ \hline \hline
pattern type          & $c_{4}$ single duty,              & $c_{4}$ single duty, $b_{3}= b_{4},$ & $c_{4}$ single duty, $b_{3}= b_{4},$  \\  
                      &  $b_{3}\ne b_{4}$                 & only 2 distinct $c$'s                & 3 distinct $c$'s     \\ \hline
$\pi$ has form        & $W_{1}\,a\,c\,W_{2}$              & $W_{1}\,a\,c\,b\,W_{2}$               & $W_{1}\,a\,c\,b\,W_{2}$       \\ \hline 
$W_{1}$ comprises     & $[b+1,n]\b \{c\}$                 & $[b-1,n]\b \{b,c\}$                  & $[b+1,n]\b \{c\}$  \\	\hline
$W_{2}$ comprises     & $[1,b]\b \{a\} $                  & $[1,b-2]\b \{a=b-3\} $               & $[1,b-1]\b \{a\} $   \\ \hline	\hline
$\rho=$               & reduce$(W_{1}\,c)$                & reduce$(W_{1}\,c)$                   & reduce$(W_{1}\,c)$  \\ \hline
ranges over           & $P^{(3)}(n-b)$                    & $\{\rho\in P^{(0)}(n-b+1):$          & $P^{(1)}(n-b)$    \\  
                      &                                   & last $\ne \min\}$                    &          \\  \hline
                      & $\binom{n-b-3}{1}2^{n-b-4}+$      &                                      & \\
a set of size         & $\binom{n-b-3}{2}2^{n-b-5}+$      & $2^{n-b-1}$                          & $\binom{n-b-2}{1}2^{n-b-3}$ \\
					  & $\binom{n-b-4}{3}2^{n-b-7}$       &                                      &   \\  \hline\hline
$\sigma=$             & reduce$(W_{2})$                   & reduce$(W_{2})$                      & reduce$(W_{2})$ \\  \hline
ranges over           & $P^{(0)}(b-1)$                    & $P^{(0)}(b-3)$                       & $P^{(0)}(b-2)$  \\  \hline       
a set of size         & $2^{b-2}$                         & $2^{b-4}$                            & $2^{b-3}$   \\  \hline\hline
sum over              & $b\in[2,n-4]$                     & $b\in[4,n-1]$                        & $b\in[3,n-3]$   \\ \hline
                      & $\binom{n-4 }{2}2^{n-6}+$         &                                      &   \\
case total            & $\binom{n-4}{3}2^{n-7}+$          & $\binom{n-4}{1}2^{n-5}$              &  $\binom{n-4}{2}2^{n-6}$      \\  
					  & $\binom{n-5}{4}2^{n-9}$           &                                      &                           \\ \hline
\end{tabular}
\\ 

\newpage 
\begin{tabular}{|c|}  \hline
$\{\pi\in P^{(4)}(n)\}$ continued \\ \hline
\end{tabular}

\vspace*{2mm}
\begin{tabular}{|c|c|c|}  \hline
 Case & 4 & 5 \\ \hline \hline
pattern type          & $c_{4}$ double duty               & $c_{4}$ double duty  \\  
                      & only 5 distinct pattern letters   & $>$ 5 distinct pattern letters \\  \hline
$\pi$ has form        & $W_{1}\,b\!-\!1\,c_{2}\,b\!-\!2\,c_{4}\,W_{2}$    & $W_{1}\,b\!-\!1\,b\!-\!2\,c\,W_{2}$ \\ \hline
$W_{1}$ comprises     & $[b+1,n]\b \{c_{2},c_{4}\}$       & $[b+1,n]\b \{c\}$        \\	\hline
$W_{2}$ comprises     & $[1,b-3]\cup \{b\} $              & $[1,b-3]\cup \{b\} $       \\ \hline	\hline
$\rho=$               & reduce$(W_{1}\,c_{2}\,c_{4})$     & reduce$(W_{1}\,c)$  \\ \hline
ranges over           & $\{\rho\in P^{(0)}(n-b)$: last 2 $\downarrow\}$ & $P^{(2)}(n-b)$     \\  \hline       
a set of size         & $2^{n-b-2}$                       & $\binom{n-b-3}{1}2^{n-b-4}+\binom{n-b-3}{2}2^{n-b-5}$ \\ \hline \hline        
$\sigma=$             & reduce$(W_{2})$                   & reduce$(W_{2})$      \\  \hline
ranges over           & $P^{(0)}(b-2)$                    & $P^{(0)}(b-2)$  \\  \hline       
a set of size         & $2^{b-3}$                         & $2^{b-3}$  \\  \hline\hline
sum over              & $b\in[3,n-2]$                     & $b\in[3,n-4]$           \\ \hline
case total            & $\binom{n-4}{1}2^{n-5}$           & $\binom{n-5}{2}2^{n-7}+\binom{n-5}{3}2^{n-8}$    \\ \hline                                   
\end{tabular}
\\ 
\vspace*{2mm}
\begin{tabular}{|c|c|c|}  \hline
 Case & 6 & 7 \\ \hline \hline
pattern type          & $c_{4}$ triple duty               & $c_{4}$ quadruple duty  \\  \hline 
$\pi$ has form        & $W_{1}\,a_{2}\,a_{3}\,a_{4}\,c\,W_{2}$    & $W_{1}\,a_{1}\,a_{2}\,a_{3}\,a_{4}\,c\,W_{2}$ \\ \hline
$W_{1}$ comprises     & $[b+1,n]\b \{c\}$                 & $[b+1,n]\b \{c\}$        \\	\hline
$W_{2}$ comprises     & $[1,b]\b \{ a_{2},a_{3},a_{4} \} $& $[1,b]\b \{a_{1}, a_{2},a_{3},a_{4} \} $  \\ \hline	\hline
$\rho=$               & reduce$(W_{1}\,c)$                & reduce$(W_{1}\,c)$  \\ \hline
ranges over           & $P^{(1)}(n-b)$                    & $P^{(0)}(n-b)$     \\  \hline       
a set of size         & $\binom{n-b-2}{1}2^{n-b-3}$       & $2^{n-b-1}$ \\ \hline \hline        
$\sigma=$             & reduce$(W_{2})$                   & reduce$(W_{2})$      \\  \hline
ranges over           & $P^{(0)}(b-3)$                    & $P^{(0)}(b-4)$  \\  \hline       
a set of size         & $2^{b-4}$                         & $2^{b-5}$  \\  \hline\hline
sum over              & $b\in[4,n-3]$                     & $b\in[5,n-1]$           \\ \hline
case total            & $\binom{n-5}{2}2^{n-7}$           & $\binom{n-5}{1}2^{n-6}$    \\ \hline                                   
\end{tabular}
\\ 
\vspace*{2mm} 
\begin{tabular}{|c|}  \hline
Grand Total:\ $2\binom{n-4}{1}2^{n-5}+3\binom{n-4}{2}2^{n-6}+\binom{n-4}{3}2^{n-7}+\binom{n-5}{3}2^{n-8}+\binom{n-5}{4}2^{n-9}$ \\ \hline
\end{tabular}
\end{center}

\end{document}